\documentclass{amsart}
\usepackage{epsfig,graphicx,color}

\begin{document}
\def\E{\ifmmode{\mathbb E}\else{$\mathbb E$}\fi} 
\def\N{\ifmmode{\mathbb N}\else{$\mathbb N$}\fi} 
\def\R{\ifmmode{\mathbb R}\else{$\mathbb R$}\fi} 
\def\Q{\ifmmode{\mathbb Q}\else{$\mathbb Q$}\fi} 
\def\C{\ifmmode{\mathbb C}\else{$\mathbb C$}\fi} 
\def\H{\ifmmode{\mathbb H}\else{$\mathbb H$}\fi} 
\def\Z{\ifmmode{\mathbb Z}\else{$\mathbb Z$}\fi} 
\def\P{\ifmmode{\mathbb P}\else{$\mathbb P$}\fi} 
\def\T{\ifmmode{\mathbb T}\else{$\mathbb T$}\fi} 
\def\SS{\ifmmode{\mathbb S}\else{$\mathbb S$}\fi} 
\def\DD{\ifmmode{\mathbb D}\else{$\mathbb D$}\fi} 

\renewcommand{\a}{\alpha}
\renewcommand{\b}{\beta}
\renewcommand{\d}{\delta}
\newcommand{\D}{\Delta}
\newcommand{\e}{\varepsilon}
\newcommand{\g}{\gamma}
\newcommand{\G}{\Gamma}
\newcommand{\la}{\lambda}
\newcommand{\La}{\Lambda}
\newcommand{\n}{\nabla}
\newcommand{\var}{\varphi}
\newcommand{\s}{\sigma}
\newcommand{\Sig}{\Sigma}
\renewcommand{\t}{\tau}
\renewcommand{\th}{\theta}
\renewcommand{\O}{\Omega}
\renewcommand{\o}{\omega}
\newcommand{\z}{\zeta}

\newcommand{\ben}{\begin{enumerate}}
\newcommand{\een}{\end{enumerate}}
\newcommand{\be}{\begin{equation}}
\newcommand{\ee}{\end{equation}}
\newcommand{\bea}{\begin{eqnarray}}
\newcommand{\eea}{\end{eqnarray}}
\newcommand{\bc}{\begin{center}}
\newcommand{\ec}{\end{center}}

\newtheorem{thm}{Theorem}[section]
\newtheorem{cor}[thm]{Corollary}
\newtheorem{lem}[thm]{Lemma}
\newtheorem{prop}[thm]{Proposition}
\newtheorem{ax}{Axiom}
\newtheorem{conj}[thm]{Conjecture}

\theoremstyle{definition}
\newtheorem{defn}{Definition}[section]

\theoremstyle{remark}
\newtheorem{rem}{\rm\bfseries{Remark}}[section]
\newtheorem*{notation}{Notation}

\newtheorem{ques}{\rm\bfseries{Question}}[section]
\newtheorem{cons}[rem]{\rm\bfseries{Construction}}
\newtheorem{exm}[rem]{\rm\bfseries{Example}}

\setcounter{page}{1}




\def\eea{\end{eqnarray*}}
\def\bea{\begin{eqnarray*}}
\def\Bbb{\bf}
\def\QQ{{\Bbb Q}}
\def\FF{{\Bbb F}}
\def\X{{\mathcal{X}}}
\def\ZZZ{{\mathcal{Z}}}
\newcommand\Oh{{\mathcal O}}
\newcommand\sC{{\mathcal C}}
\newcommand\sA{{\mathcal A}}
\newcommand\sD{{\mathcal D}}
\newcommand\sF{{\mathcal F}}
\newcommand\sG{{\mathcal G}}
\newcommand\sI{{\mathcal I}}
\newcommand\sJ{{\mathcal J}}
\newcommand\sL{{\mathcal L}}
\newcommand\sB{{\mathcal B}}
\newcommand\sN{{\mathcal N}}
\newcommand\sK{{\mathcal K}}
\newcommand\sM{{\mathcal M}}
\def\de{{\delta}}
\def\De{{\Delta}}
\def\ga{{\gamma}}
\def\a{{\alpha}}
\def\be{{\beta}}
\def\Ga{{\Gamma}}
\def\PPP{{\Bbb P}}
\def\ZZ{{\Bbb Z}}
\def\NN{{\Bbb N}}

\newcommand{\Proof}{{\it Proof. }}
\newcommand{\QED}{{\hfill $Q.E.D.$}}

\newtheorem{teo}[thm]{Theorem}
\newtheorem{df}[thm]{Definition}
\newtheorem{ex}[thm]{Example}
\newtheorem{oss}[thm]{Remark}

\newcommand{\LL}{\ensuremath{\mathbb{L}}}
\newcommand{\F}{\ensuremath{\mathbb{F}}}
\newcommand{\M}{\ensuremath{\mathbb{M}}}
\newcommand{\hol}{\ensuremath{\mathcal{O}}}
\newcommand{\HH}{\ensuremath{\mathbb{H}}}
\newcommand{\PP}{\ensuremath{\mathbb{P}}}
\newcommand{\RRR}{\ensuremath{\mathcal{R}}}
\newcommand{\BB}{\ensuremath{\mathcal{B}}}
\newcommand{\FFF}{\ensuremath{\mathcal{F}}}
\newcommand{\A}{\ensuremath{\mathcal{A}}}
\renewcommand{\D}{\ensuremath{\mathcal{D}}}
\newcommand{\CCC}{\ensuremath{\mathcal{C}}}
\newcommand{\EE}{\ensuremath{\mathcal{E}}}
\newcommand{\HHH}{\ensuremath{\mathcal{H}}}
\newcommand{\V}{\ensuremath{\mathcal{V}}}
\newcommand{\I}{\ensuremath{\mathcal{I}}}
\newcommand{\SSS}{\ensuremath{\mathcal{S}}}
\newcommand{\VV}{\ensuremath{\mathbb{V}}}
\newcommand{\MM}{\ensuremath{\mathfrak{M}}}
\newcommand{\NNN}{\ensuremath{\mathfrak{N}}}
\newcommand{\gS}{\ensuremath{\mathfrak{S}}}

\newcommand{\CC}{\ensuremath{\mathbb{C}}}
\newcommand{\RR}{\ensuremath{\mathbb{R}}}
\newcommand{\MMM}{\ensuremath{\mathcal{M}}}
\newcommand{\BBB}{\ensuremath{\mathcal{B}}}
\newcommand{\VVV}{\ensuremath{\mathcal{V}}}
\newcommand{\ra}{\ensuremath{\rightarrow}}
\newcommand{\rla}{\leftrightarrow}
\newcommand{\AAA}{\ensuremath{\mathcal{A}}}

%
%
\newcommand{\blau}[1]{{\color{blue} #1}}
\newcommand{\rot}[1]{{\color{red} #1}}
\newtheorem{lemma}[thm]{Lemma}

\newcommand{\ON}{\operatorname}
\newcommand{\inv}{^{^{-1}}}
\newcommand{\invv}{^{-2}}
\newcommand{\br}{\ON{Br}}

\newcommand{\cg}{\gamma}
\newcommand{\Cg}{\Gamma}
\newcommand{\eps}{\epsilon}
\newcommand{\varth}{\vartheta}
\renewcommand{\k}{\kappa}
\renewcommand{\l}{\lambda}
\newcommand{\oo}{\omega}

\newcommand{\cutoff}[1]{}

\newcommand{\id}{\ON{id}}

\allowdisplaybreaks

%
%



\title[Moduli spaces and monodromy invariants.]
{Moduli spaces of surfaces and monodromy invariants.}

\author[Catanese \it et al.]{Fabrizio Catanese, Michael L\"onne, Bronislaw Wajnryb}

\thanks{}

\address{Prof. Fabrizio Catanese\\
Lehrstuhl Mathematik VIII\\
Universit\"at Bayreuth, NWII\\
           D-95440 Bayreuth, Germany}
\email{Fabrizio.Catanese@uni-bayreuth.de}

\address{Prof. Michael L\"onne\\
    Mathematisches Institut der\\
Universit\"at G\"ottingen,\\
Bunsenstrasse 3\\
           D-37073 G\"ottingen, Germany
}
\email{Michael.Loenne@uni-bayreuth.de}

\address{Prof. Dov Bronislaw Wajnryb \\
Department of Mathematics and Applied Physics\\
Technical University of Rzeszow\\
     Rzeszow, Poland
}
\email{dwajnryb@prz.edu.pl}

\begin{abstract}
Bidouble covers $\pi : S \ra Q: = \PP^1 \times \PP^1$ of the quadric
and their natural deformations are
parametrized by connected families $\mathfrak N_{a,b,c,d}$ depending
on four positive integers
$a,b,c,d$. We shall call these surfaces abcd-surfaces.
In the special case where $ b=d$ we call them abc-surfaces.

Such a Galois covering $\pi$  admits a small  perturbation
yielding a general 4-tuple covering of $Q$ with branch curve $\De$,
and two natural Lefschetz fibrations obtained from a small perturbation
of the composition $ p_i \circ \pi, \ i=1,2$, $p_i$ being the i-th
projection of $Q$ onto $\PP^1 $.

The first and third author
showed that the respective mapping class group factorizations corresponding to
the first  Lefschetz fibration are equivalent for two
abc-surfaces with the same values of $ a+c, b$, a result which implies the
diffeomorphism of two such surfaces.

We report on  a more general result of the three authors implying
           that the first braid monodromy factorization corresponding to $\De$
determines  the three integers $a,b,c$ in the case of abc-surfaces.
We provide in this article a new proof of the non equivalence of two
such factorizations
for different values of $a,b,c$.

We finally show that, under certain conditions, although the first
    Lefschetz fibrations are equivalent for two
abc-surfaces with the same values of $ a+c, b$,
the second Lefschetz fibrations need not be equivalent.

These results rally around the question whether
abc-surfaces with fixed values of $ a+c, b$,
although diffeomorphic but not deformation equivalent, might  be not
canonically symplectomorphic.

\end{abstract}


\maketitle

\tableofcontents

\section{Introduction: Moduli spaces of surfaces of general type:
    deformation, differentiable and
symplectic types.}

Let  $X$ be a smooth complex projective
variety, i.e.,   $X$ is a smooth
compact complex submanifold of $\PP^n : =  \PP^n_{\CC}$  of complex
dimension $d$,
and  $X$ is connected.

Observe that $X$ is a compact oriented real manifold of real
dimension $2d$, so, for instance, a complex surface gives
rise to a real 4-manifold.

For $d=1$, $X$ is a complex  algebraic curve (a real surface, called
also Riemann surface) and its basic invariant
is the genus $g= g(X)$, defined as
the dimension of the vector space $H^0 (\Omega^1_X)$ of holomorphic 1-forms.

The situation for curves is `easy', since the genus determines the
topological and the
differentiable manifold underlying $X$: its intuitive meaning is the
`number of handles'
that one must add to a sphere in order to obtain $X$ as a topological space.

The rough classification of curves is the following :

\begin{itemize}

\item $g=0$ : $X \cong \PP^1_{\CC}$, topologically $X$ is a sphere $S^2$
of  real dimension $2$.
\item $g=1$ : $X \cong \CC / \Gamma$ , with $\Gamma$ a discrete subgroup
$\cong \Z^2$: $X$ is called an elliptic curve, and topologically we have
a real 2-torus $ S^1 \times S^1$.
\item $ g\geq 2 $ : then we have a `curve of general type', and topologically
we have a connected sum of $g$ 2-tori.
\end{itemize}

Every curve is a deformation of a special hyperelliptic curve, the
non singular projective model
of the affine curve
$$ \{ (x,z) \in \CC^2 | z^2 = (x-1) (x-2) \dots ( x - (2g +2)) \}. $$

The main theorem for curves says more precisely that we have a {\bf
Moduli space} $\frak
M_g$ which parametrizes the isomorphism classes $[C]$ of compact
complex curves $C$
of genus $g$.

$\frak M_g$ is a Zariski open set of a complex projective
variety, is connected, and has  complex dimension
$ (3g-3) + a(g)$, where $a(g)$ is the complex dimension
of the group of complex automorphisms ($ a(0)= 3, a(1) = 1, a(g) = 0 
,\  {\it for} \
g \geq 2$).

When we pass to complex dimension $d=2$, some features do generalize,
others do not.

The first generalizations of the genus were given by Clebsch, Noether,
Enriques and
Castelnuovo through the dimensions of certain vector spaces of
holomorphic tensor differential forms
$$q(X): = dim_{\CC} H^0 (\Omega^1_X),  p_g(X): = dim_{\CC}H^0 (\Omega^d_X),
P_m(X): = dim_{\CC}H^0
({\Omega^d_X}^{\otimes m}).$$
$q(X)$ is called the irregularity of $X$, $p_g(X)= P_1(X)$ is called
the geometric genus,
while $P_m(X)$ is called the m-th plurigenus.

These invariants suffice to give the Castelnuovo-Enriques classification of
algebraic surfaces.

The rough classification of projective surfaces $S$ is the following :

\begin{itemize}

\item $P_{12} = 0$ : $S$ is birational to a product $ C \times \PP^1_{\CC}$,
where $C$ is a complex curve of genus $q(S)$.
\item $P_{12} = 1$ : the surface is birational to a surface $S$
which is an analogue
of an elliptic curve, in the sense that there is a holomorphic section of $H^0
({\Omega^2_S}^{\otimes 12})$ which is nowhere vanishing.
\item $P_{12} \geq 2$ and $S$ is properly elliptic, i.e.,  the
rational map $\phi_{12}$ associated to the sections of $H^0
({\Omega^2_S}^{\otimes 12})$ maps to a curve and the general fibres
are elliptic curves.
\item $P_{12} \geq 2$ and $S$ is of general type, i.e., $\phi_{12}$
maps birationally
to a surface $X$, called the canonical model of $S$.
\end{itemize}

Surfaces of special type are well understood, but surfaces of general
type still offer
us a lot of intriguing and fascinating open problems.
In each birational class there is for them a unique smooth projective surface
(up to isomorphism) $S$, called minimal model, with the property that
every birational
holomorphic map  $S$ onto another smooth  surface $S'$ is a biholomorphism.
This minimal model offers therefore also a unique topological and
differentiable model.
As does the genus of an algebraic curve, the  above numerical invariants are determined
by the topological structure of $S$:
\begin{itemize}

\item  $2 q(S)$ is twice the first Betti number of $S$, $2 q(S) = b_1(S)$,  and
\item
if one defines  $\chi(S) : =
1 - q(S) + p_g(S)$, we have Noether's formula $ 12 \chi(S)= K^2_S + e
(S)$, where $e(S) =
c_2(S)$ is the topological Euler Poincar\'e characteristic, while
$K^2_S = c_1^2(S) $; moreover
\item
the  signature $\sigma(S) : = b^+(S)  - b^-(S)$ of
the intersection form $q_S : H^2(S, \Z) \ra \Z$ equals, by virtue of Hodge's
index theorem, 
$$\sigma(S) =  1/3 (K_S^2 - 2 e(S))$$ and finally we have Kodaira's
Riemann-Roch formula (holding  for surfaces of general type) which says that
\item
the plurigenus $P_m(S)$, for $ m \geq 2$, equals $ \frac{ m (m-1)}{2}
K_S^2 + \chi(S)$.
\end{itemize}

Things  seem `on the surface' to generalize nicely: because    also for
algebraic surfaces
of general type there exist similar moduli spaces $\frak M_{x,y}$
(by the results of \cite{bom}, \cite{gieseker}).

Here, $\frak M_{x,y}$ is quasi projective and parametrizes
isomorphism classes of minimal
(smooth projective) surfaces
of general type $S$ such that $ \chi (S) = x, K^2_S = y$.

The fact that $\frak M_{x,y}$ is quasi-projective implies that
$\frak M_{x,y}$ has a finite number of connected components, which parametrize
deformation classes
of surfaces of general type; and, by a classical theorem of Ehresmann
(\cite{ehre}),
deformation equivalent varieties are diffeomorphic.

Hence, fixed the two numerical invariants $\chi (S) = x, K^2_S = y$,
which are determined by the topology of $S$,
we have a finite number of
differentiable types, and a fortiori a finite number of topological types.

These two numbers  $ \chi (S), K^2_S $ are determined by the
topology. But they do not determine the topology, since the fundamental group
is not encoded in such invariants (and is not even invariant under
the action of the
absolute Galois group, cf. \cite{serre}, and  \cite{bcg} for new examples).

For this reason one likes to restrict to the case of simply connected surfaces.
In this case, these two numbers almost determine the topology.
This follows by
Michael  Freedman's big Theorem of 1982 (\cite{free}), showing that there are
indeed at most  two topological structures if  the  surface $S$ is
assumed to be simply connected.

Topologically, our 4-manifold is then obtained  from  very simple
building blocks,
one of them being the K3 surface, where:

\begin{df}
A {\bf K3} surface is a smooth surface of degree 4 in $\PP^3_{\CC}$.
\end{df}

Now, a complex manifold carries a natural orientation
corresponding to the complex structure, and, in general, given an oriented
differentiable manifold $M$, $M^{opp}$ denotes the same manifold, but endowed
with the opposite orientation. This said, one can explain the corollary
of Freedman's theorem for the topological manifolds underlying
simply connected (compact) complex surfaces as follows:

\begin{itemize}
\item
     if $S$ is EVEN, i.e., the intersection form on $H^2 (S, \ZZ)$ is even:
then $S$ is homeomorphic to a connected sum of copies of
      $ \PP ^1_{\CC} \times \PP ^1_{\CC}$ and copies of  a $K3$ surface
if the signature $\sigma(S)$ is negative, and of copies of
$ \PP ^1_{\CC} \times \PP ^1_{\CC}$ and of copies of a $K3$ surface with
reversed orientation if the signature is positive.

\item
if $S$ is ODD : then $S$ is a connected sum of copies of
$\PP ^2_{\CC} $ and of ${\PP ^2_{\CC}}^{opp} .$
\end{itemize}

Kodaira and Spencer defined quite generally (\cite{k-s58}) for  compact complex
manifolds
$X, X'$  the equivalence relation called {\bf deformation equivalence}: this,
      for  surfaces of general type, means that the corresponding
isomorphism classes
      yield points in the same connected component
of the moduli space $\frak M$.

The cited theorem of Ehresmann guarantees that DEF $\Rightarrow $ DIFF:
indeed a bit more holds, namely,
deformation equivalence implies the existence of a diffeomorphism carrying
the canonical class $K_X$ to the canonical class $K_{X'}$.

In the 80' s  the work of Simon Donaldson (\cite{don1},
\cite{don2},\cite{don3},\cite{don4}, see also \cite{d-k}) showed that
homeomorphism and diffeomorphism differ drastically
for projective surfaces.

The introduction of  Seiberg-Witten invariants
(see \cite{witten},\cite{don5},\cite{mor})
showed then more easily that  a diffeomorphism $\phi : S \ra S'$
between minimal surfaces of general type satisfies
$\phi^* (K_{S'}) = \pm  K_S$.

Contrary to the  conjecture made by the first author in \cite{katata},
Friedman and Morgan made the

FRIEDMAN-MORGAN'S SPECULATION (1987) (see \cite{f-m1}):

   differentiable  and deformation equivalence coincide for surfaces.
   
We abbreviate this speculation with the acronym

  DEF =  DIFF.

The question was answered negatively in every
possible way (\cite{man4}, \cite{k-k}, \cite{cat4}, \cite{c-w}, \cite{bcg},
see also \cite{catcime} for a  rather comprehensive survey):

\begin{teo}
(Manetti '98, Kharlamov -Kulikov 2001, Catanese 2001, Catanese - Wajnryb 2004,
Bauer- Catanese - Grunewald 2005 )

The Friedman- Morgan speculation
does not hold true.
\end{teo}

The counterexamples obtained by Catanese and  Wajnryb are
the only ones which are simply connected, and one  also has
(as also in Manetti's examples)
non deformation equivalent surfaces $S, S'$ such that there exists
    an orientation preserving  diffeomorphism $\phi : S \ra S'$ with
$\phi^* (K_{S'}) =   K_S$.

The  simply connected examples  used, called {\bf abc-surfaces}, are
a special case of a class of surfaces which the first author introduced in 1982
(\cite{cat1}), namely, bidouble covers of the quadric
and their natural deformations.

We shall briefly describe these surfaces in the next section.

But, before, let us explain why we believe that these surfaces deserve
further attention. Recall the

\begin{df}
A pair $(X,\omega)$ of a real manifold $X$, and of a real differential
2- form $\omega$ is called a {\bf Symplectic pair} if

i) $\omega$ is a symplectic form, i.e., $\omega$ is closed ( $ d \omega  = 0$)
       and  $\omega$ is nondegenerate at each point (thus $X$ has even
real dimension).

A symplectic pair $(X,\omega)$ is said to be {\bf integral} iff the De Rham
cohomology class of $\omega$ comes from $ H^2 ( X, \Z)$, or, equivalently,
there is a complex line bundle $L$ on $X$ such that $\omega$ is a first Chern
form of $L$.

ii) An  almost complex structure $J$ on $X$ (a differentiable endomorphism
of the real tangent bundle of $X$ satisfying $ J^2 = - 1$) is said to be
compatible with $\omega$  if
$ \omega (Jv, Jw) = \omega (v, w)$, and
the quadratic form
$\omega ( v, J v)$ is strictly positive definite.

\end{df}

For long time (before the  examples
of Kodaira and Thurston) the basic examples of
symplectic
manifolds were given by symplectic submanifolds
    of K\"ahler manifolds,
in particular  of the flat space $\CC^n$ and of projective
space $\PP^N_{\CC}$,
endowed with the Fubini-Study form
$$\omega_{FS} := \frac{i}{2 \pi } \partial
\overline{\partial} log |z|^2.$$

\bigskip

It was observed recently by the first author (\cite{cat02}, \cite{cat06}):

\begin{teo}
      A minimal surface of general type $S$ has a symplectic structure
$(S,\omega)$,
unique up to symplectomorphism, and invariant for smooth deformation,
with class$(\omega) = K_S = -  c_1(S).$

This symplectic structure is called the
{\bf  canonical symplectic structure}.
\end{teo}

The above result is, in the case where $K_S$ is ample, a
rather direct consequence of the well known Moser's
lemma (\cite{moser}), since then it suffices
to pull-back $1/m$ of the Fubini-Study metric by an
    embedding of $S$ by the sections of $H^0 (S, {\Omega^2_S}^{\otimes m})$.

It is the case where $K_S$ is not ample which needs a non trivial  proof.

An important consequence of the techniques of symplectic
approximation of singularites
employed in these papers is
the following theorem (\cite{cat02},
\cite{cat09}, see also \cite{catcime} for a survey including the
basics concerning the
construction of the Manetti surfaces)

      \begin{teo}
Manetti's surfaces  yield examples of surfaces of general type which are
not deformation equivalent but are canonically symplectomorphic.
\end{teo}

The following questions are still wide open.

{\bf Questions:} 1) Are there (minimal) surfaces of general type which
are orientedly
       diffeomorphic through a diffeomorphism carrying the canonical class
to the canonical
class, but, endowed with their canonical symplectic structure,
are not
canonically   symplectomorphic?

2) Are there such   simply connected examples ?

3) Are there diffeomorphic abc-surfaces which are canonically symplectomorphic,
thus yielding a  counterexample to Can. Sympl = Def
also in the simply connected case?
Are there diffeomorphic abc-surfaces which are not canonically
symplectomorphic,
    thus providing an answer to question 1) ?

The difficult problem
of understanding the canonical symplectic structures of abc-surfaces
(\cite{clw}) is relying on the work of Auroux and Katzarkov (\cite{a-k1}),
who described some invariants of integral symplectic pairs of real
dimension 4. We postpone to a later section the description of their method.

\section{abcd- Surfaces and their moduli spaces}

All the curves are deformations of hyperelliptic curves $C$ , and the
nature of the double cover
$ \pi : C \ra \PP^1$ allows an easy understanding of the topology of $C$.

Double covers of the projective plane $\PP^2$ or of $\PP^1 \times \PP^1$ are however too simple
surfaces, with very special behaviour; while it turned out that   bidouble covers of  $\PP^1 \times \PP^1$, 
introduced  in \cite{cat1}, exhibited a quite non special behaviour 
concerning properties of the  moduli spaces of surfaces of general type.

Bidouble covers of the quadric are smooth projective complex surfaces $S$
endowed with a (finite) Galois
covering $ \pi : S \ra Q : = \PP^1 \times \PP^1$ with Galois group
$(\Z /2 \Z)^2$.

More concretely, they are defined by a single pair of equations

\begin{eqnarray*}
      & z^2= f_{(2a,2b)}(x_0 , x_1; y_0 , y_1)\\
      & w^2=g_{(2c,2d)}(x_0 , x_1; y_0 , y_1)\\
\end{eqnarray*}

\noindent where $ a,b,c,d\in\N^{\geq3}$ and the notation $f_{(2a,2b)}$ denotes that
$f$ is a bihomogeneous polynomial, homogeneous of degree $2a$ in the
variables $x$,
respectively  of degree $2b$ in the variables $y$.

These surfaces are simply connected and minimal of general type,
and they  were introduced in \cite{cat1} in order to show that the
      moduli spaces
$\frak M_{\chi,K^2}$ of smooth minimal surfaces of general type $S$ with
$K^2_S = K^2 , \chi(S): = \chi (\hol_S) = \chi$
need not be equidimensional or irreducible (and indeed the same holds
for the open and closed subsets $\frak M_{\chi,K^2}^{00}$ corresponding to
simply connected surfaces).

In fact, for the dimension of the moduli space  $\frak M_{\chi,K^2}$ one can only give lower and upper
bounds, as shown by the followig  theorem (cf. \cite{cat1})

\begin{thm}
The dimension of an irreducible component $\frak M'_{\chi,K^2}$ of $\frak M_{\chi,K^2}$ is subject to the following inequalities:
$$10 \chi - 2 K^2 \leq  dim (\frak M'_{\chi,K^2} ) \leq    10 \chi  +
3  K^2  + 108.$$
\end{thm}

And the deformations of  bidouble covers of the quadric were used to show that
the moduli spaces $\frak M_{\chi,K^2}^{00}$ could have as many 
irreducible components of different dimensions as possible.

These irreducible components were determined as follows:
given  our four integers $ a,b,c,d\in\N^{\geq3}$,  consider
the so called natural deformations of these bidouble covers
of {\em simple type}  $(2a,2b),(2c,2d)$,
defined by equations \footnote{in the following formula,
a polynomial of negative degree is identically zero. }

$$  (***) \ \ \  \ z_{a,b}^2 =  f_{2a,2b} (x,y)  + w_{c,d} \phi_{2a-c, 2b-d}(x,y)$$
       $$ w_{c,d}^2 = g_{2c,2d}(x,y) +  z_{a,b}  \psi_{2c-a,2d -b}(x,y)$$

where f, g, $\phi, \psi$, are bihomogeneous polynomials ,  belonging to
        respective vector spaces of sections of line bundles (recall that we defined $Q : = {\PP^1\times \PP^1}$):

        $ f \in H^0({Q}, {\hol}_{Q}(2a,2b)) $,

$ \phi \in  H^0({Q}, {\hol}_{Q}(2a-c, 2b-d))$ and

$ g \in  H^0({Q}, {\hol}_{Q}(2c,2d))$,

$\psi \in  H^0({Q}, {\hol}_{Q}(2c-a, 2d - b))$.

In this way one defines an open\footnote{deformation theory shows
that this set is open}
   subset $\NNN'_{a,b,c,d}$ of the moduli space,
whose closure $\overline{\NNN'}_{a,b,c,d}$
is an irreducible component of $\frak M_{\chi,K^2}$, where
      $ \chi = 1 + (a-1)(b-1) + (c-1)(d-1) + (a+c-1)(b+d-1)$, and $K^2 =
8 (a+c-2)(b+d-2)$.

We record here for further use:
$ \chi  - \frac{1}{8} K^2 =  ab + cd.  $

These calculations can be derived from the formula 
$\hol_S ( K_S) = \pi^* \hol_Q (a+c-2, b+d-2)$, which was later used in
\cite{cat3} to show
that the index of divisibility $r(S)$ of $K_S$ in $H^2(S, \Z)$ equals the
$ G.C.D. (a+c-2, b+d-2) $. This result  showed
that many of these irreducible components  $\overline{\NNN'}_{a,b,c,d}$
would belong
to distinct connected components of the moduli space
(since the divisibility index $r(S)$ is constant on a connected component).

The  {\bf abc-surfaces} are obtained as the special case where $ b = d$,
and the upshot of \cite{c-w} is that, once the  values of the integers  $b$
and  $a+c$ are fixed, one obtains diffeomorphic
surfaces.

It is relevant to observe that the family $\NNN'_{a,b,c,d}$ is contained
as a dense open subset in a larger irreducible family $\NNN_{a,b,c,d}$,
where  coverings of degree 4 of other rational surfaces
$\F_{2h} : = Proj (\hol_{\PP^1} \oplus \hol_{\PP^1} (2h))$
have also to be considered.

\bigskip

In \cite{cat1} it was shown that two irreducible components of the
moduli space,
such as  $\overline{\NNN'}_{a,b,c,d}$,
even when they have the same dimension, are distinct as soon as the 4-tuples
are not trivially equivalent by an obvious $(\Z/2)^2$-symmetry (exchanging the two factors of ${\PP^1\times \PP^1}$,
or the two branch divisors). In fact an irreducible component 
defines a surface over the generic point,
and if two families coincide, then they should define isomorphic surfaces over
the generic point. In particular, these two surfaces should be isomorphic 
as surfaces over the algebraic closure of the function field of the generic point.

In this particular case, the corresponding surfaces were shown to be non isomorphic
since this isomorphism would yield an isomorphism of the canonical images,
and the canonical images were shown to have non equivalent singularities.

In general, it is much more difficult to show that two irreducible components yield distinct connected
components of the moduli space.

In practice, in order to show  that a certain open set $\mathfrak N$
of the moduli space is a connected component, one  has to prove that the set 
$\mathfrak N$ is also closed.

In order to do so one has to look at all families of surfaces parametrized by a smooth curve
$T$, and such that all the surfaces
$S_t , \ t \in T \setminus \{ t_0\}$ belong to $\mathfrak N$: if one can 
   show that also the limit
surface $S_{ t_0}$ stays in $\mathfrak N$ then  $\mathfrak N$ is then closed,
hence a connected component.

Of course, if two surfaces $S, S'$ are not diffeomorphic (e.g., if their divisibiity indices of the canonical class satisfy $r(S)
\neq r(S')$),
then they lie in different connected components.

Boris Moishezon proposed instead to use, as an invariant of a 
connected component $\mathfrak N$
of the moduli space, the geometry of pluricanonical projections, 
which we shall describe
in the final section.

\section {Diffeomorphisms of abc-surfaces and Lefschetz pencils}

The   main new result of \cite{c-w} was the following

\begin{teo}\label{diffabc} {\bf (Catanese -Wajnryb )}
Let $S$ be an $(a,b,c)$ - surface (i.e., an $(a,b,c,b)$ surface) and $S'$ be an
$(a+1,b,c-1)$-surface. Moreover, assume that $a,b,c-1 \geq 2$. Then
$S$ and $S'$ are diffeomorphic.
\end{teo}

This result  was then put together  with the following more technical result,
in order to exhibit simply connected surfaces which are diffeomorphic but not
deformation equivalent.

\begin{teo} \label{nondef}  {\bf (Catanese -Wajnryb )}
Let  $S$, $S'$ be simple bidouble
covers of  ${\PP}^1
\times  {\PP}^1  $ of respective
        types ((2a, 2b),(2c,2b), and (2a + 2k, 2b),(2c - 2k,2b), and assume

\begin{itemize}
\item
        (I) $a,b,c, k$
are strictly  positive even integers with $ a, b, c-k \geq 4$
\item
(II)  $ a \geq 2c + 1$,
\item
(III) $ b \geq c+2$  and either
\item
(IV1) $ b \geq 2a + 2k -1$ or\\
(IV2) $a \geq b + 2$
\end{itemize}
Then  $S$ and  $S'$ are not deformation equivalent.
\end{teo}

       The  second theorem uses techniques which have been developed
in a  series of papers by the  first author and by Marco Manetti during
a long period of time, and it belongs more to algebraic geometry
than to geometry and topology.

The essential thrust is to show, under the above technical conditions,
that $\NNN_{a,b,c,b}$ is indeed a connected component of
the moduli space.

The substantial inequality is the inequality (II), which guarantees that
for all natural deformations one has $\psi \equiv 0$ (since a
polynomial of negative degree
is identically zero), i.e., the natural deformations have the following form

$$  (***) \ \ \  \ z_{a,b}^2 =  f_{2a,2b} (x,y)  + w_{c,b} \phi_{2a-c,b}(x,y)$$
       $$ w_{c,b}^2 = g_{2c,2b}(x,y) .$$

The above equations show that the covering $\pi : S \ra Q$ is an
iterated double covering,
i.e., we first take an intermediate  double covering $Y$ determined
by taking the square root
$w$ of $g$, and then we obtain (through the quadratic equation for
$z$) $S$ as a double covering
$S \ra Y$.

The main idea is that, while taking limits,  we can recover a smooth
rational surface
$\F_{2h}$ as the  limit of $Q = Y/ (\Z/2)= (S/ (\Z/2))/ (\Z/2)$.

We mention here however that the following question is rather open in general:

{\bf Question:} 1) For which values of ${a,b,c,d}$ is $\NNN_{a,b,c,d}$
a connected component of the moduli space?

A more detailed expositions for both theorems can be found
in the   Lecture Notes of the C.I.M.E. courses
`Algebraic surfaces and symplectic 4-manifolds' (see especially
\cite{catcime}).

The above result in differential topology \ref{diffabc}  is based instead on a
refinement of Lefschetz theory obtained by Kas (\cite{kas}).

This refinement allows us to encode the differential topology of a 4-manifold
$X$ Lefschetz fibred over $\PP^1_{\CC}$ (i.e.,
$f : X \ra \PP^1_{\CC}$ has the property that all the fibres are smooth and
connected, except for a finite number which have exactly one nodal singularity)
into an equivalence class of a factorization of the identity
in  the Mapping class group $\MMM ap_g$ of a compact curve
$C_0$ of genus $g$.

Recall that the mapping class group, introduced by Max Dehn(\cite{dehn})
in the 30's, is defined for each orientable manifold $M$
as
$$ \MMM ap (M) : =  Diff ^+(M) / Diff^0 (M), $$
where $Diff^+ (M)$ is the group of orientation preserving diffeomorphisms,
while $Diff^0 (M)$ is the connected component of the identity,
the so called subgroup of the diffeomorphisms which are
{\bf isotopic to the identity}.

In the case of a complex curve of genus $g$, Dehn found that the group
was generated by the so called Dehn twists around simple  closed
loops
   $\Ga$,
which are diffeomorphisms equal to the identity outside of a
neighbourhood of $\Ga$,
and  rotate $\Ga$ by the antipodal map.

In our situation, let $F : X \ra \PP^1_{\CC}$ be a Lefschetz fibration, and
let $b_1, \dots b_m$ be the critical values of $f$.
For each $b_i$ the fibre $C_i : = f^{-1} (b_i)$ has a nodal singularity $P_i$,
and the local monodromy around $b_i$ is given, by a theorem of Picard
and Lefschetz,
by the Dehn twist around the vanishing cycle in a nearby smooth fibre
$C_{b'_i}$.

If one fixes a base point $b_0$ in the base $\PP^1_{\CC}$, and simple
paths $\ga'_i$ not crossing each other
from $b_0 $ to $b'_i$, one obtains diffeomorphisms, well defined up
to isotopy, of the standard fibre $C_0 : = C_{b_0}$
with the  smooth fibres $C_{b'_i}$, and then a sequence of Dehn
twists $\tau_1, \dots \tau_m$.

Assuming that the paths $\ga'_i$ are ordered in counterclockwise
order, by adding a small circle
around the point $b_i$ run counterclockwise starting from $b'_i$, one
obtains a geometric
basis  $\ga_i$ for the fundamental group
$$\pi (\PP^1_{\CC} \setminus \{  b_1, \dots b_m \}, b_0) = \langle
\ga_1, \dots \ga_m | \ga_1 \cdot  \dots  \cdot \ga_m = 1\rangle.$$

By virtue of the relation $ \ga_1 \cdot  \dots  \cdot \ga_m = 1$ we
obtain a sequence of Dehn twists
 $ \tau_1,   \dots  \tau_m $ on $C_0= C_{b_0}$  whose product $ \tau_1\cdot
\dots  \cdot  \tau_m $
equals the identity. We obtain a factorization of the identity in
$\MMM ap_g$ once we fix a diffeomorphism
of the fibre $C_0 = C_{b_0}$  with a fixed compact complex curve of genus $g$.

This means that this factorization is only defined up to
simultaneous conjugation in the group
$\MMM ap_g$, i.e., one can replace, for $\phi  \in \MMM ap_g$, each $\tau_i$
by $$ \tau_i^{\phi} : = \phi^{-1} \tau_i \phi .$$

Moreover, there is another ambiguity for the factorization, due
namely to the choice of the
generators $\ga_1, \dots \ga_m$ for the fundamental group: this
ambiguity amounts to an action of the braid group in $m$ strings on
the factorization,
leading to the so called Hurwitz equivalence for factorizations, and
which will be described in the next section.

Said this, one can state the theorem of Kas:
\begin{thm}
Two Lefschetz fibrations $F_1, F_2$ of the same genus $g \geq 2 $ are
diffeomorphic (there is a diffeomorphism $\varphi$ with $F_1 = F_2
\circ \varphi$)
iff the corresponding factorizations in the mapping class group $\MMM
ap_g$ are equivalent via Hurwitz equivalence and simultaneous
conjugation.

Moreover, if, for $j=1,2$, $S_j$ is a complex surface and $F_j : S_j 
\ra \PP^1$ is homotopic
to a holomorphic map, then the diffeomorphism $\varphi$ sends the 
canonical class
$ - c_1(S_1) \in H^2 (S_1, \ZZ)$ to the canonical class $ - c_1(S_2)$.
\end{thm}

\proof
We have only added a
second part to the statement, based on the Seiberg Witten theorem
mentioned above, and on the fact that $\varphi$ sends fibres to 
fibres in an orientation preserving fashion.
Hence $\varphi$  is orientation preserving, sends the class of the 
fibre to the class of a fibre,
$ c_1(S_1)$ to $ \pm c_1(S_2)$: we conclude that the sign is $+1$, 
since the intersection number of a fibre
with $ - c_1(S_j)$ equals $ 2g - 2 > 0$.

\QED

The concrete application of Kas' theorem present two big difficulties:

1) how to write such a monodromy factorization?

2) how to verify that two factorizations are equivalent ?

3) how to verify that two factorizations are not equivalent ?

\medskip

For the first question,we shall show how we can calculate such
factorizations, using
the reality principle: trying to write fibrations which are defined
over $\R$, and possibly with all the critical
values being real (if not, they come into conjugate pairs).
   In this case the paths we choose can be explicitly described.
In particular, it helps to write the branch curves as smoothings of
rather singular but explicitly given
curves, union of graphs of rational functions.

The second question is quite hard, but, at least in the case of
abc-surfaces, a big help came from a
simple geometric idea.

In fact, consider our surfaces as bidouble covers $Z$

$$  (***) \ \ \  \ z_{a,b}^2 =  f'_{2a,2b} (x,y)  $$
       $$ w_{c,b}^2 = g'_{2c,2b}(x,y) $$

branched on a union of horizontal and of vertical lines: in this case
the smooth surface
$S$ is the blow up of $Z$ at the nodes lying over the nodes of $C' =
\{ f'(x,y) =0 \}$
and of $D' = \{ g'(x,y) =0 \} $.

Changing the bidegree of $f'$ from $(2a, 2b)$ to $(2a + 2, 2b)$, and
the bidegree of $g'$ from $(2c, 2b)$ to $(2c - 2, 2b)$ amounts geometrically
to removing two vertical lines from the curve $D'$ and putting them back
into the curve $C'$. Here we exploit the fortunate circumstance that $b=d$:
we take the horizontal lines of $C'$ as the lines $ y=j, \ j=1, \dots 2b$,
and the horizontal lines of $D'$ as the lines $ y=- j, \ j=1, \dots  2b$.

Consider then the manifold with boundary $M$, resolution of the
branched cover of
$ \Delta \times \PP^1$, where $\Delta$ is the unit disk in $\C$,
described by
$$   \ \ \  \ z^2 =  \Pi_1^{2b} ( y - i) $$
       $$ w^2 = (x^2 - \epsilon^2) \Pi_1^{2b} ( y + i).$$

Then $S'$ is obtained from $S$ removing $M$ and attaching it back on
the boundary of
$\Delta$ after acting with the automorphism $\psi$ of order two such that
$ \psi (x,y,z,w) =  (x, -y, w,z)$.

The result is that the factorizations for $S$ and $S'$ differ just  in this
beginning part, where
the factors have been conjugated by $\psi$.
Geometrically, one says that the fibrations of $S$ and $S'$ are
obtained as two different fibre sums:
i.e., one adds to
the same factorization (the final part)  a beginning
factorization, once as it is, and another time
conjugated by $\psi$.
   There is a very nice lemma, due to Auroux (mentioned also
parenthetically by Kas), which gives a
criterion for equivalence of these two beginning pieces of the
respective factorizations:

\begin{lem}\label{auroux} {\bf ( Auroux)} Let $\tau$ be a Dehn twist and let
$F $  be a factorization of a central element
$ \phi \in \MM ap_g$,
$ \tau_1 \circ \tau_2 \circ \dots \circ \tau_m = \phi$.

If there is a factorization $F'$ such that $F$ is Hurwitz equivalent
to $\tau \circ F'$, then $ (F)^{\tau} : = \tau^{-1} F \tau$ is Hurwitz equivalent
to $F$.

In particular, if $F$ is a factorization of the identity,
$\psi = \Pi _h \tau'_h$, and $\forall h \exists F'_h$ such that
$F $ is Hurwitz equivalent to $ \tau'_h \circ F'_h$, then the fibre sum with the Lefschetz
pencil associated with $F$ yields the same Lefschetz pencil
as the fibre sum twisted by $\psi$.

\end{lem}
\Proof

If $\cong$ denotes Hurwitz equivalence, then
$$ (F)^{\tau} \cong \tau \circ (F')^{\tau} \cong F' \circ \tau
\cong (\tau)^{(F')^{-1}} \circ F' =  \tau \circ F' \cong F.$$

\QED

\begin{cor} \label{cor}
Notation as above, assume that $F$,
$ \tau_1 \circ \tau_2 \circ \dots \circ \tau_m = \phi$,
is a factorization of the Identity
and that $\psi$ is a product of the Dehn twists $\tau_i$ appearing
in $F$. Then a fibre sum with the Lefschetz
pencil associated with $F$ yields the same result
as the same fibre sum twisted by $\psi$.
\end{cor}

\Proof
We need only to verify that for each $h$ , there is $F'_h$ such that
$F \cong \tau_h \circ F'_h$.

But this is immediately obtained by applying $h-1$ Hurwitz moves,
the first one between $\tau_{h-1}$ and $\tau_h$, and proceeding further
to the left till we obtain $\tau_h$ as first factor.

\QED

Hence in \cite{c-w}  the problem was reduced to proving that the
isotopy class of $\psi$ was a product
of certain Dehn twists appearing in the first part
of the factorization.

Verifying isotopy of these diffeomorphisms  was accomplished by constructing chains of loops in the
complex fibre curve $C_0$, which lead to a dissection of $C_0$ into open cells.
The problem was then solved  choosing several associate Coxeter elements
to express the given diffeomorphism $\psi$ as a product
of those Dehn twists.

In the next section we shall instead concentrate on question 3),
which is related to
our recent work.

\section{How to distinguish factorizations up to Hurwitz equivalence,
and stable equivalence}

In this section the set up will be more algebraically oriented: the
reason for doing so
stems from the fact that in our applications we shall deal with two
different type of groups, not only mapping
class groups, but also braid groups, and moreover we shall have to
deal, in the second case, with two
different types of equivalence relations, the first one being the
equivalence relation generated by
simultaneous conjugation and by Hurwitz equivalence, the second one,
introduced by Auroux and Katzarkov, is called stable equivalence and
is more general.

We shall look  carefully into the problem of distinguishing
classes of factorizations, and for this reason we shall
introduce a  technical
novelty allowing a new effective method for disproving
stable-equivalence.

Assume that we have a group $G$.

Then a factorization in $G$ of length $m$ is nothing else
than an element in the Cartesian product $G^m$.

The product function associates to such an element
$ (\a_1,    \dots  \a_m ) \in G^m$
its product $ \a_1 \cdot  \dots  \cdot \a_m = a$.

The inverse image of the element $a$ under the product function is
called  the set of {\bf  factorizations of $a$ in $G$ of length $m$}
and we write such a  factorization of $a$ in $G$ with the notation
$$(*) \  \a_1 \circ   \dots  \circ \a_m = a.$$

There is a natural action of the group $ Aut (G)$ on the set of  length $m$
factorizations, in particular $G$ acts by conjugation, hence
the Centralizer $Z(a)$ acts on the factorizations of $a$ by
simultaneous conjugation,
$$  (\a_1 \circ   \dots  \circ \a_m = a) \mapsto   (\a_1^g \circ
\dots  \circ \a_m^g = a), \
\forall g \in Z(a), $$
where as before we set $a^g : = g^{-1} a g$.

Hurwitz equivalence of factorizations is the equivalence generated by
Hurwitz moves,
where the i-th Hurwitz move leaves all the $\a_j$'s with $j \neq i,
i+1$ unchanged,
while it replaces $\a_{i+1}$ with $\a_{i}$, and $\a_{i}$ with
$\a_{i}\a_{i+1}\a_{i}^{-1}$,
so that the product of the factorization remain unchanged.

Hurwitz equivalence classes are nothing else than orbits for the action of
the braid group $\BB_m$ on the set of factorizations of a given element $a$.

Recall (cf.  e.g. \cite{birman}):
\begin{teo}
The braid group $\BB_n$ is the subgroup of
$\MM ap ( \C - \{1, \dots n\})$ given by the  diffeomorphisms which are the
identity outside of a circle with  centre the origin
and radius $2n$.
It has the Artin presentation
$$ < \sigma_1,  \dots \sigma_{n-1} |
   \sigma_i  \sigma_j =  \sigma_j
   \sigma_i  \ {\rm for }\  |i-j| \geq 2, \ \sigma_i  \sigma_{i+1}  \sigma_i =
\sigma_{i+1}  \sigma_i \sigma_{i+1} \ ,  \forall i \leq n-2>. $$
\end{teo}

Through this isomorphism the standard Artin generators
$\sigma_j$ of  $\BB_n$ ( $ j=1, \dots n-1$) correspond to  standard
half-twists: $\sigma_j$ is the diffeomorphism of $(\C, \{1, \dots n\} 
)$ isotopic to the
homeomorphism  of $(\C, \{1, \dots n\} )$
which is a rotation of 180 degrees on the disc with centre $j+1/2$
and radius $1/2$ ,
on a circle with the same  centre and radius $ 1/2 + t/4$ is
the identity for $ t \geq 1$,
and a  rotation
of angle $180 (1-t)$ degrees if  $0 \leq  t \leq 1$.

Since $\BB_n$ is a subgroup of
$\MM ap ( \C - \{1, \dots n\})$ we have an obvious action of
$\BB_n$ on the free group $\pi_1( \C - \{1, \dots n\})$.
Taking as base point  the complex number $p := -2ni$, one has a geometric basis
$\gamma_1, \dots \gamma_n$ on which $\BB_n$  acts through the

{\bf  Hurwitz action of the braid group}

\begin{itemize}
\item
$\sigma_i (\gamma_i) = \gamma_{i+1}$
\item
$\sigma_i (\gamma_i \gamma_{i+1}) = \gamma_i \gamma_{i+1}$  hence
$\sigma_i (\gamma_{i+1}) = \gamma_{i+1} ^{-1}\gamma_i \gamma_{i+1}$
\item
$\sigma_i (\gamma_j ) = \gamma_j $ for $ j \neq i, i+1$.
\end{itemize}

The Hurwitz action leaves the product
$\gamma_1  \gamma_2 \dots    \gamma_n $ invariant.

We can summarize the above discussion in the following

\begin{df}
Let $G$ be a group, and considered the partition of its  Cartesian product
$G^n$ given by the product map, which to an ordered  n-uple $(g_1,
g_2, \dots g_n)$
assigns the
product $g := g_1\cdot  g_2 \dots \cdot  g_n \in G$.

The subsets of the partition of $G^n$
are called the  {\bf factorizations of an element
$ g \in G$}.
The orbits of the action of  $\BB_n$, which preserves the partition of $G^n$,
are called {\bf Hurwitz equivalence classes of  factorizations}.

Two factorizations are said to be {\bf M-equivalent} (this means: 
equivalent by Hurwitz
equivalence
and simultaneous conjugation, and the M- stands for monodromy) if 
they are in the same  $\BB_n \times G$-orbit,
where $G$ acts by (we use the notation previously introduced $ a^b :
= b^{-1} a b$)
the transformation

   $$ Int_b \colon (g_1 \circ g_2  \circ \dots  \circ g_n   \mapsto
   (g_1)^b \circ (g_2)^b  \circ \dots  \circ (g_n)^b.$$
\end{df}

\begin{oss}
Given a group homomorphism $\phi : G \ra G'$,and a factorization in 
$G$, $ g_1 \circ g_2  \circ
\dots  \circ g_n  = a$, we  obtain a corresponding factorization 
$$\phi(g_1) \circ \phi(g_2)  \circ
\dots  \circ \phi(g_n )) = \phi (a),$$ through the application of $\phi^n$.

$\phi^n$ preserves Hurwitz equivalence, and simultaneous conjugation. 
Therefore, a basic method
to disprove equivalence of factorizations is to disprove equivalence 
of homomorphic factorizations
in a simpler group $G'$.

This is essentially the main underlying idea.
\end{oss}

There are some obvious invariants for a factorization class:

\begin{prop}\label{inv}
Consider the set of Hurwitz equivalence classes of factorizations
of a fixed element $ a \in G$:
$$(*) \  \a_1 \circ   \dots  \circ \a_m = a.$$

Then the basic invariants of a Hurwitz  factorization class are:

i)  the subgroup $H$ generated by the
elements $\a_1,   \dots   \a_m$

ii) the function, defined on the set $\sC _H$ of $H$-conjugacy classes in $H$,
assigning to each conjugacy class $A \in \sC_H$ the integer $\sigma (A)\in \N$
counting how many of the elements $\a_1,   \dots   \a_m$
belong to $A$.

Consider instead the set of M-equivalence classes of factorizations
of a fixed element $ a \in G$:
$$(*) \  \a_1 \circ   \dots  \circ \a_m = a.$$

Then the basic invariants of an M- factorization class are,
$Z(a)$ being the centralizer of the element $ a \in G$:

1) the $Z(a)$-conjugacy class of the subgroup $H$ generated by the
elements $\a_1,   \dots   \a_m$, in particular the isomorphism class
of $H$ as an abstract group.

2) the $Z(a)$ class of the function $\sigma \colon \sC_H \ra \N$,
in particular, the cardinalities $\sigma^{-1}(n)$, for $n \in \N$.

\end{prop}

In the holomorphic case the above invariants are well defined; however,
in the symplectic world Auroux and Katzarkov  had to introduce
a broader equivalence of factorizations, where the length of the
factorization is not fixed a priori. They called it $m$-equivalence,
but we prefer to call it stable-equivalence, or S-equivalence.

\begin{df}
Consider a group $G$, an element $a \in G$,
and a set $\BB$ of elements $\be_j \in G$.

We  define {\bf stable}-equivalence with respect to $\BB$ simply by
considering the equivalence relation generated by Hurwitz equivalence,
simultaneous conjugation and creation/cancellation of consecutive factors
$  \be_j \circ \be_j^{-1}$.

\end{df}

Then  {\bf refined  invariants } of the stable equivalence class of
a factorization are obtained
as follows:

I) let $H$ be the subgroup of $G$ generated by the $\a_i$'s
and let $\hat{H}$ be the subgroup of $G$ generated by the $\a_i$'s
and by the $\be_j$'s (in our case, $H$ will be called the monodromy group,
and $\hat{H}$ the stabilized monodromy group).

Then the $Z(a)$-conjugacy class of $\hat{H}$ is a first invariant,
in particular the isomorphism class of $\hat{H}$.
\smallskip

II) Let $\sC_H$ be the set of $H$-conjugacy classes $A$ in the  group $H$,
and let $\hat{\sC}$ be the set of $\hat{H}$-conjugacy classes in the  group
$\hat{H}$, so that we have a natural map $\sC \ra \hat{\sC}$,
with $A \mapsto \hat{A}$.

Write $\hat{\sC}$ as a disjoint union $\hat{\sC}^+ \cup \hat{\sC}^-
\cup \hat{\sC}^0$,
where $\hat{\sC}^0$ is the set of conjugacy classes of elements $a$ which
are conjugate
to their inverse $ a^{-1}$, and where $\hat{\sC}^-$ is the set of the inverse
conjugacy classes of the classes in $\hat{\sC}^+$.

Associate now to the factorization $ \a_1 \circ   \dots  \circ \a_m = 1$
the function $ s : \hat{\sC}^+ \ra \Z$ such that $ s(c)$, for $ c \in
\hat{\sC}^+$,
is the algebraic number of occurrences of $c$ in the sequence of
conjugacy classes
of the $ \a_j $'s (i.e., an occurrence of $c^{-1}$ counts as $-1$ for $ s(c)$).

The $Z(a)$ class of the function $ s : \hat{\sC}^+ \ra \Z$ is  our 
second and most
important invariant.  From it one can derive an easier invariant,
i.e., the $Z(a)$ class of the function induced by the absolute value 
$|s|$ on the set
  $\hat{\sC}^{+-}$ of pairs of distinct conjugacy classes, one inverse
  of the other. In particular, the values $|s|^{-1}(n)$ for $n \in \N$.

\begin{oss}
The calculation of the function $ s : \hat{\sC}^+ \ra \Z$
        presupposes however a detailed knowledge of the group $  \hat{H}$.
For this reason in our paper \cite{clw} we resorted to using  some coarser
derived invariant.
\end{oss}

{\bf Substrategy I.}

Assume that  we can write $\{\a_1, \dots, \a_m\} $ as a disjoint
union $\sA_1 \cup \sA_2 \cup \sD \cup \sA'_1 \cup \sA'_2 $, such that
the set $\sA_j, \ j=1,2,$
is contained in a conjugacy class $A_j \subset H$,
the set $\sA'_j, \ j=1,2,$ is contained in the conjugacy class
$A_j^{-1} \subset H$.
Assume that the elements  in $\sA_1 \cup \sA_2$ are contained in a conjugacy
class $C$ in $G$ such that $ C \cap C^{-1} = \emptyset$ but that
the set $\sD$ is disjoint from the union of the two conjugacy classes of $G$,
$ C \cup C^{-1} $.

If we then prove that $\hat{A_1} \neq \hat{A_2} $ (this of course implies
$A_1 \neq A_2 $) we may assume w.l.o.g.  that
$\hat{A_1} , \hat{A_2} \in \hat{\sC}^+$, and  then the unordered pair
of positive numbers $( |s(\hat{A_1} )|, |s(\hat{A_2}| )$  is our derived
numerical invariant of the factorization (and can easily be calculated from the
cardinalities of the four sets  as $(||\sA_1| - |\sA'_1|| , ||\sA_2|
- |\sA'_2||)$.

\bigskip
{\bf Substrategy II.}

This is the strategy to show that $\hat{A_1} \neq \hat{A_2} $
and goes as follows.

Assume further that we have another  subgroup of $G$, $ \tilde{H} \supset \hat{H}$,
and a group homomorphism $\rho :  \tilde{H}  \ra \Sigma$ such that
$\rho ( \BB ) = 1$.  Assume also that the following key property holds.

{\bf Key property.}

        For each element $ \a_j \in A_j \subset H$ there exists
an element $ \tilde{\a}_j \in \tilde{H}$
        such that $$ \a_j = \tilde{\a}_j^2 ,$$
and moreover that this element is unique  in $\tilde{H}$ (a fortiori, it
will suffice that it  is unique in $G$).

Proving that  $ \rho (\tilde{\a}_1)$ is not conjugate to
$ \rho (\tilde{\a}_2)$ under the action of $ \rho (H) = \rho (\hat{H})
\subset \Sigma$
shows finally that $\a_1$ is not conjugate to $\a_2$ in $\hat{H}$.

Since, if there is $ h \in \hat{H}$ such that $\a_1 = h^{-1} \a_2 h$,
then $\tilde{\a}_1^2  = \a_1 = (h^{-1} \tilde{\a}_2 h)^2$, whence
$\tilde{\a}_1  = h^{-1} \tilde{\a}_2 h$, a contradiction.

In our concrete case, we are able to determine the braid monodromy group $H$,
and we observe that, since we have a so called cuspidal factorization,
all the factors $\a_i$ belong to only four conjugacy classes
in the group $G$ (the classes of $\sigma_1$,$\sigma_1^2$,
$\sigma_1^{-2}$,$\sigma_1^3$
in the braid group), corresponding geometrically to vertical tangencies of the
branch curve, respectively positive nodes, negative nodes and cusps.
Three of these classes are positive and only one is negative (the one
of $\sigma_1^{-2}$).

We see in \cite{clw} that  the nodes belong to  conjugacy classes $A_1, A_2,
A_1^{-1}, A_2^{-1}
\subset H$, the positive nodes belonging to  $A_1 \cup  A_2$,
and we show, using a representation $\rho$ of  a certain subgroup $\tilde{H}$
of `liftable' braids (i.e., braids which centralize the monodromy
homomorphism, whence are liftable to the mapping class group of
a curve $D_0$ associated to the monodromy homomorphism) into a symplectic group
$\Sigma$ with $\Z / 2$ coefficients,
that the classes $\hat{A_1}, \hat{A_2}$ are distinct in $\hat{H}$. Moreover,
these are classes in  $\sC^+$, since these are positive classes in the
braid group.
We calculate then easily the above function for these two conjugacy classes,
i.e., the pair of numbers $ ( s(\hat{A_1}), s(\hat{A_2})) $.

\section{Bidouble covers of the quadric and their symplectic perturbations.}

In order to see in more detail how the algebraic considerations of 
the previous sections apply to our
situation,
let us go back to our equations

$$  (***) \ \ \  \ z_{a,b}^2 =  f_{2a,2b} (x,y)  + w_{c,d} \Phi_{2a-c, 2b-d}(x,y)$$
       $$ w_{c,d}^2 = g_{2c,2d}(x,y) +  z_{a,b}  \Psi_{2c-a, 2d-b}(x,y)$$

where $f,g$ are bihomogeneous polynomials as before, and instead we shall
allow $\Phi,
\Psi$, in the case where for instance the degree relative to $x$ is negative,
to be antiholomorphic.  In other words, we allow $\Phi,
\Psi$ to be sections of certain line bundles which are dianalytic
(holomorphic or antiholomorphic) in each variable $x, y$.

We can more generally consider the direct sum $\VV$ of two complex line bundles
$\LL_1 \oplus \LL_2$
on a compact complex manifold  $X$, and the subset  $Z$
of $\VV$ defined by the following pair of equations,

\begin{eqnarray*}
           & z^2= f(x) + w \
\Phi (x)\\
           & w^2=g(x) + z \
\Psi (x)\\
\end{eqnarray*}

where $f,g$ are respective holomorphic sections of the line bundles
$\LL_1^{\otimes 2},  \LL_2^{\otimes 2}$ (we shall also write
$ f \in H^0 (\hol_X(2 L_1)),  g \in H^0 (\hol_X(2 L_2))$, denoting by
$L_1, L_2$ the associated
Cartier divisors), and where  $\Phi$ is a
differentiable section
of  $\LL_1^{\otimes 2} \otimes \LL_2^{\otimes -1}$,
$\Psi$ is a differentiable section
of  $\LL_2^{\otimes 2} \otimes \LL_1^{\otimes -1}$.

The key observation is that, for $\Phi \equiv \Psi \equiv 0$, then
$\pi : Z \ra X$ is a Galois cover with group $(\Z/2\Z)^2$,
while, if $\Phi ,  \Psi $ are not $\equiv  0$, then we obtain a 
covering whose monodromy group
is the symmetric group $\gS_4$.

We have the following Lemmas from \cite{clw}.

\begin{lem}\label{smallperturbation}
Assume that the two divisors $\{ f=0\}$ and $\{ g=0\}$ are smooth and
intersect transversally.
Then, for $|\Phi| < < 1$, $|\Psi| < < 1$, $Z$ is a smooth submanifold
of $\VV$, and the projection
$\pi : \VV \ra X$ induces  a finite covering  $ Z \ra X$ of degree
$4$, with ramification divisor (i.e., critical set) $ R : = \{ 4 zw =
\Phi \Psi \}$, and
with branch divisor (i.e., set of critical values) $\De = \{ \de(x) = 0\}$,
where $$ - \frac{1}{16} \de = - f^2 g^2 - \frac{9}{8} fg (\Phi
\Psi)^2 + (\Psi)^2 f^3 +
          (\Phi)^2 g^3 + \frac{27}{16^2} (\Phi \Psi)^4.$$
\end{lem}

We want to analyse now the singularities that $\De$ has, for a
general choice of
$\Phi, \Psi$. The singularities which are easily spotted deserve a 
special definition.

\begin{df}
A point $x$ is said to be a {\bf trivial singularity} of $\De$ if it
is either a point where $ f = \Phi = 0$, or a point where $ g = \Psi = 0$.

$Z$ is said to be {\bf mildly general } if

(*)  at any point where $ f = \Phi = 0$,
we have $ g \cdot\Psi \neq 0$, and symmetrically at any point where
$ g = \Psi = 0$, we have $ f \cdot\Phi \neq 0$
\end{df}

Besides the trivial  singularities, we show in general that the other 
singularities are located in a neighbourhood of the
points with$ f(x) = g(x)
=0$.

\begin{lem}
If $Z$ is mildly general, and  $|\Phi| < < 1$, $|\Psi| < < 1$,
the singular points of $\De$ which are not
the trivial singularities of $\De$
occur only for points  $p \in Z_x$ ($Z_x$ being the fibre over the point $x$)
of multiplicity $3$ which lie over  arbitrarily  small neighbourhoods
of the points $x'$ with $ f(x') = g(x')
=0$.
\end{lem}

We can describe more precisely the singularities of $\De$ in the case where
$X$ is a complex surface.

\begin{prop}
\label{singularities}
Assume that $X$ is a compact complex surface, that the two divisors
$\{ f=0\}$ and $\{ g=0\}$ are
smooth and intersect transversally in a set $M$ of $m$ points.
Then, for $|\Phi| < < 1$, $|\Psi| < < 1$, and for $Z$  mildly general,
$Z$ is a smooth submanifold of
$\VV$, and the projection
$\pi : \VV \ra X$ induces  a finite covering  $ Z \ra X$ of degree
$4$, with smooth orientable ramification divisor $ R : = \{ 4 zw =
\Phi \Psi \}$, and
with branch divisor  $\De $ having as singularities precisely

1) $ 3m$ cusps lying
(in triples) in an arbitrarily small
neighbourhood of $M$, and moreover

2) the trivial singularities, which are nodes if
the curve $\{ f=0\}$
intersects transversally $\{ \Phi = 0\}$, respectively
if the curve  $\{ g=0\}$
intersects transversally $\{ \Psi = 0\}$.
\end{prop}

Since $R$ is orientable, it follows that if the trivial singularities
of $\De$ are nodes,
then the intersection number of the two branches is $\pm 1$; the local
selfintersection number is  exactly equal to $+1$ when the two orientations of
the two branches combine to yield the natural (complex) orientation of $X$.

We can calculate the precise number of singularities, and of the
positive, respectively negative nodes of $\De$, in the special case we
are interested in,
namely, of a perturbed bidouble cover of the quadric.

\begin{df}\label{separately}
Given  four integers $ a,b,c,d\in\N_{\geq3}$,
the so called {\bf dianalytic perturbations  of simple bidouble covers}  are
the {\rm 4}-manifolds
defined by equations
\begin{eqnarray*}
           & z^2= f_{(2a,2b)}(x_0 , x_1; y_0 , y_1) + w \
\Phi _{(2a-c,2b-d)}(x_0 , x_1; y_0 , y_1)\\
           & w^2=g_{(2c,2d)}(x_0 , x_1; y_0 , y_1) + z \
\Psi _{(2c-a,2d-b)}(x_0 , x_1; y_0 , y_1)\\
\end{eqnarray*}
where $f,g$ are bihomogeneous  polynomials of respective bidegrees
$(2a,2b)$,

          $(2c,2d)$,
          and where  $\Phi, \Psi$ are
polynomials in the  ring $$\C[ x_0 , x_1, y_0 , y_1, \bar{x_0 },
\bar{x_1},\bar{y_0 },\bar{y_1}],$$
          bihomogeneous   of
respective bidegrees
$(2a-c,2b-d)$,  $(2c-a,2d-b)$ (the ring is bigraded here by setting the
degree of  $\bar{x_i }$ equal to
$(-1,0)$ and the degree of  $\bar{y_i }$ equal to
$(0, -1)$ ).

We choose moreover $\Phi, \Psi$ to belong to the subspace where
all monomials are either separately holomorphic or antiholomorphic,
i.e., they admit no factor of the form $x_i \bar{x_j }$ or
of the form $y_i \bar{y_j }$.
\end{df}

\begin{oss}\label{antihol}
Indeed it is convenient to  pick $\Phi$ (respectively : $\Psi$) to be a product
$$\Phi = \Phi_1(x)  \Phi_2 (y),$$ where $ \Phi_1$ is a product
of
linear forms, either all holomorphic or all antiholomorphic,
and similarly for  $\Phi_2$ (resp. : for $\Psi_1$, $\Psi_2$).

\end{oss}

As explained in \cite{catcime}(pages 134-135) to an antiholomorphic
homogeneous polynomial
$ P(\bar{x_0 },
\bar{x_1})$ of degree $m$ we associate a differentiable section $p$ of
the tensor power $\LL^{\otimes m}$ of the tautological (negative) subbundle,
such that $ p (x) : = P( 1, \bar{x})$ inside a big disc $B(0,r)$
in the complex line
having centre at the origin and  radius $r$.

We make the following assumptions on the polynomials $f,g, \Phi, \Psi$:

\begin{enumerate}\label{hyp}
\item
The algebraic curves $C : = \{ f = 0 \}$ and $D : = \{ g = 0 \}$ are smooth
\item
$C$ and $D$ intersect transversally at a finite set $M$ contained in
the product $B(0,r)^2$  of two big discs in $\C$, and at
these points both curves have non vertical tangents
\item
for both curves $C$ and $D$ the first projection on $\PP^1$ is a
simple covering
and moreover the vertical tangents of $C$ and $D$ are all distinct
\item
the associated perturbation is sufficiently small and mildly general
\item
the trivial singularities are nodes with non vertical tangents,
contained in $B(0,r)^2$
\item
the cusps of $\De$ have non vertical tangent.
\end{enumerate}

Considering then  the branch curve $\De$ of the perturbed bidouble cover,
the defining equation $\de$ is bihomogeneous of bidegree $(4
(a+c), 4 (b+d))$.

We have that $\De$ has exactly $ k : = 12 (ad + bc)$ cusps, coming from the
          $ m : = 4 (ad + bc)$ points of the set $ M = C \cap D$.

          If $\Phi, \Psi$
are holomorphic, then $\De$ has exactly
$$\nu : = 4 ( 2 ab + 2 cd - ad - bc)$$
(positive) nodes.

In general, the above  number $\nu$ equals
the difference $\nu^+ - \nu^-$ between the number of positive and
the number of negative nodes.

In fact the nodes of $\De$ occur only for the trivial singularities.

A trivial singularity $ f = \Phi = 0$ yields a node with  the same tangent cone
    as $$ f^2 - \Phi^2 \cdot g(x_0).$$

By remark \ref{antihol} we may assume that there are local holomorphic
coordinates $f, \hat{\Phi} $ such that either $\Phi = \mbox{ unit
}\cdot \hat{\Phi}$
or $\Phi = \mbox{ unit }\cdot \overline {\hat{\Phi}}$.

In the first case we get the tangent cone of a holomorphic node, hence
a {\bf positive} node, in the second case we
get the tangent cone of an antiholomorphic node, hence
a {\bf negative} node.

We summarize the results in the following

\begin{lem}
The number $t_f$ of vertical tangents for the curve $C$ is
$  t_f = 4 ( 2 ab - a) $, while the  number $t_g$ of vertical
tangents for the curve $D$ is
$  t_g = 4 ( 2 cd - c) $. For a general choice of $\Phi, \Psi$, the
number of vertical tangents
$t$ of $\De$ equals $ t = 2 t_f + 2 t_g + m$, where $ m : = 4 (ad + bc)$.
The genus of the Riemann surface $R$ equals $ g(R) =  1 + 16 (a+c)(b+d)
- 4 (a+b+c+d) - k -\nu $,
where $ k : =  3 m = 12 (ad + bc)$ is the number of cusps of $\De$ and
$\nu : = \nu^+ - \nu^-$ is the number of nodes of $\De$, counted with sign.
\end{lem}

\begin{oss}
1) Observe that the fact that $p : R \ra\PP^1_{\C}$ has a finite
number of critical points implies
immediately that $p$ is finite, since there is then a finite set in
$\PP^1_{\C}$
such that over its complement we have a covering space. Thus $p$ is
orientation preserving
and each non critical point contributes positively to the degree of
the map $p$. Hence some of  the above calculations are obtained  from
Hurwitz' ramification formula.

2) Consider now a real analytic 1-parameter family $Z(\eta)$ such that
the perturbation terms
$\Phi_{\eta}, \Psi_{\eta} \ra 0$. Then we have a family of ramification points
$P_i({\eta})$ which tend
to the ramification points of $p_0 : R_0 \ra\PP^1_{\C}$. For $R_0$, a
double cover of
$C \cup D$, we have $2t_f$ ramification points lying over the $t_f$
vertical tangents
of $C$,  $2t_g$ ramification points lying over the $t_g$ vertical tangents
of $D$, whereas the $m$ nodes of $R_0$ are limits of $4$ ramification points,
three cusps and a simple vertical tangent, cf. the analysis made in 
\cite{c-w2}. This is the geometric reason why
$t =  2 t_f + 2 t_g + m$; it also tells us where to look for the
vertical tangents of $\De$.
\end{oss}

In practice, in order to effectively compute the braid monodromies,
it is convenient to view the  curves $C$ and $D$ as
small real deformations
of reducible real curves with only nodes (with real tangents) as singularities,
and to  consider real polynomials $\Phi, \Psi,$ such that the
trivial singularities are also given by real points. By the results of
\cite{c-w2}, from the perturbation of each proper node of $C \cup D$,
i.e., of a point of $ M = C \cap D$, we shall obtain  a real vertical  tangent,
a real cusp, and two immaginary cusps.

\newcommand{\nodalC}{C^\times}
\newcommand{\nodalD}{D^\times}
\newcommand{\nodalf}{f^\times}
\newcommand{\nodalg}{g^\times}

We begin with real polynomials $\nodalf$ and $\nodalg$
defining nodal curves $\nodalC$, respectively $\nodalD$.

For the sake of simplicity we replace here bihomogeneous
polynomials $ f(x_0, x_1 ; y_0, y_1 )$ by their restrictions to the
affine open set
$x_0 =1 , y_0 =1$ and write $ f(x,y)$ for $ f(1, x ; 1, y )$.

\unitlength=1pt
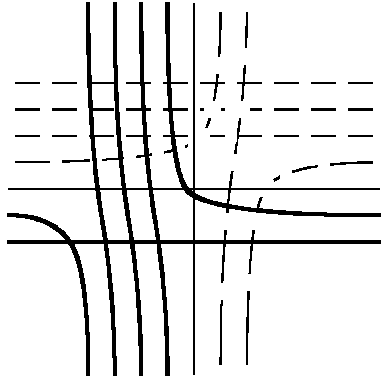
\begin{figure}
\centering
\begin{picture}(150,150)(-75,-75)
\color{white}

{\color{black}
\linethickness{.4pt}
\qbezier(-70,20)(0,20)(70,20)
\qbezier(-70,30)(0,30)(70,30)
\qbezier(-70,40)(0,40)(70,40)

\qbezier(15,10)(20,40)(20,70)
\qbezier(15,10)(10,-20)(10,-70)

\qbezier(3,20)(10,30)(10,70)
\qbezier(3,20)(-4,10)(-70,10)
\qbezier(27,0)(20,-10)(20,-70)
\qbezier(27,0)(34,10)(70,10)

{\color{white}
\linethickness{4pt}
\qbezier[10](-70,20)(0,20)(70,20)
\qbezier[10](-70,30)(0,30)(70,30)
\qbezier[10](-70,40)(0,40)(70,40)
\linethickness{4pt}
\qbezier[7](3,20)(-4,10)(-70,10)
\qbezier[3](27,0)(34,10)(70,10)

\linethickness{6pt}
\qbezier[3](27,0)(20,-10)(20,-70)
\qbezier[3](15,10)(20,40)(20,70)
\qbezier[3](15,10)(10,-20)(10,-70)

\qbezier[3](3,20)(10,30)(10,70)
}


\linethickness{.2pt}

\qbezier(-70,0)(0,0)(70,0)
\qbezier(0,-70)(0,0)(0,70)


\linethickness{1.2pt}
\qbezier(-70,-20)(0,-20)(70,-20)

\qbezier(-15,-10)(-10,-40)(-10,-70)
\qbezier(-15,-10)(-20,20)(-20,70)
\qbezier(-25,-10)(-20,-40)(-20,-70)
\qbezier(-25,-10)(-30,20)(-30,70)
\qbezier(-35,-10)(-30,-40)(-30,-70)
\qbezier(-35,-10)(-40,20)(-40,70)

\qbezier(-47,-20)(-40,-30)(-40,-70)
\qbezier(-47,-20)(-54,-10)(-70,-10)
\qbezier(-3,0)(-10,10)(-10,70)
\qbezier(-3,0)(4,-10)(70,-10)
}
\end{picture}
       \caption{Real part of $\Delta_f$ and $\Delta_g$ for
         $(a,b,c,d)=(1,2,2,1)$}
\label{delta f and g}
\end{figure}

We let then $\nodalf : =f_1\cdots f_{2b}$ defining  $f_i : = y-2i$
       for i=$2, \dots , 2b$ but setting
$$
f_1: =(y-2)\prod_{i=1}^{2a}(x-2i)+\prod_{i=1}^{2a-1}(x-2i-1),
$$
and we define similarly $\nodalg : =g_1\cdots g_{2d}$ setting $g_i=y+2i$  for
i=$2, \dots , 2d$ except that we set
$$
g_1: =(y+2)\prod_{i=1}^{2c}(x+2i)-\prod_{i=1}^{2c-1}(x+2i+1).
$$

\bigskip
\begin{oss}
\label{halfplane}
Note that the equations of $f_1$ and $g_1$ are chosen in such a way
that their zero sets are graphs of rational functions $\tilde{f}_1$,
respectively
       $\tilde{g}_1$, of
$x$ which, regarded as maps from
$\C$ to
$\C$, preserve the real line. Moreover, $\tilde{f}_1$ preserves both the upper
and lower halfplane, while $\tilde{g}_1$ exchanges them.

In fact, $\tilde{f}_1$ has no critical points on $\PP^1_{\R}$ and
$$ deg (\tilde{f}_1 |_{ \PP^1_{\R} }) \colon \PP^1_{\R }\ra \PP^1_{\R}  =
2a = deg (\tilde{f}_1 ).$$
Similarly for $\tilde{g}_1$.
\end{oss}

\unitlength=1pt
\begin{figure}
\centering
\begin{picture}(170,170)(-75,-75)
{
{
\qbezier(-70,60)(5,60)(80,60)
\qbezier(-70,70)(5,70)(80,70)
\qbezier(-70,80)(5,80)(80,80)
\color{white}
\linethickness{3pt}
\qbezier[5](-48,60)(5,60)(60,60)
\qbezier[5](-48,70)(5,70)(60,70)
\qbezier[5](-48,80)(5,80)(60,80)
}
{
\qbezier[60](45,-70)(45,10)(45,90)
\qbezier[60](55,-70)(55,10)(55,90)
\qbezier[60](65,-70)(65,10)(65,90)
}
}

{
\linethickness{.1pt}
\qbezier(-70,0)(0,0)(80,0)
\qbezier(0,-70)(0,0)(0,90)
}

\linethickness{.5pt}
\qbezier(-70,20)(0,20)(80,20)
\qbezier(-70,30)(0,30)(80,30)
\qbezier(-70,40)(0,40)(80,40)

\qbezier(15,10)(20,40)(20,90)
\qbezier(15,10)(10,-20)(10,-70)

\qbezier(3,20)(10,30)(10,90)
\qbezier(3,20)(-4,10)(-70,10)
\qbezier(27,0)(20,-10)(20,-70)
\qbezier(27,0)(34,10)(80,10)

\bezier{400}(-70,-20)(0,-20)(80,-20)

\bezier{200}(-15,-10)(-10,-40)(-10,-70)
\bezier{200}(-15,-10)(-20,20)(-20,90)
\bezier{200}(-25,-10)(-20,-40)(-20,-70)
\bezier{200}(-25,-10)(-30,20)(-30,90)
\bezier{200}(-35,-10)(-30,-40)(-30,-70)
\bezier{200}(-35,-10)(-40,20)(-40,90)

\bezier{200}(-47,-20)(-40,-30)(-40,-70)
\bezier{200}(-47,-20)(-54,-10)(-70,-10)
\bezier{200}(-3,0)(-10,10)(-10,90)
\bezier{200}(-3,0)(4,-10)(80,-10)

\end{picture} \caption{Zero sets of $\Phi$ (dashed) and $\Psi$
(dotted) added}
\label{phipsi}
\end{figure}
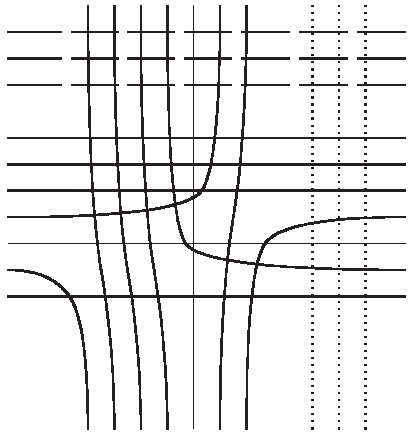

Given $\nodalf$, consider the polynomial $f : \nodalf + c_f$,
where $c_f$ is a small constant; likewise    consider  $g : \nodalg + c_g$.
Adding these small  constants to the
respective equations of
$\nodalC$ and
$\nodalD$
we get polynomials $f$ and $g$ which define
smooth curves $C$ and $D$. We have more precisely:

\begin{prop}
\label{hypprop}
If the constants $c_f$, $c_g$, $\eta$
are chosen sufficiently small, the polynomials $f, g, \eta \Phi, \eta 
\Psi$ thus obtained
satisfy  the following list of hypotheses,
cf. \pageref{hyp},  for $r >8(a+b+c+d)$:
\begin{enumerate}\label{hyp2}
\item
The algebraic curves $C = \{ f = 0 \}$ and $D = \{ g = 0 \}$
are smooth
\item
$C$ and $D$ intersect transversally at a finite set $M$ contained in
the product $B(0,r)^2$  and at
these points both curves have non vertical tangents
\item
for both curves $C$ and $D$ the first projection on $\PP^1$ is a
simple covering
and moreover the vertical tangents of $C$ and $D$ are all distinct
\item
the associated perturbation is mildly general
\item
the trivial singularities of $\Delta$
are nodes with non vertical tangents,
contained in $B(0,r)^2$
\item
the cusps of $\De$ have non vertical tangent
\end{enumerate}
\end{prop}

\section{Bidouble covers of the quadric and their  Braid Monodromy 
factorization types.}

In the previous section we showed how  to obtain a symplectic 
4-manifold which (we call a
dianalytic perturbation and) is
canonically symplectomorphic to the bidouble cover $S$ we started with.

Now  we have  that, for general choice of $f,g, \Phi,
\Psi$, the projection $\pi$ onto $Q = \PP^1
\times
\PP^1$ is {\bf generic} and its branch curve $\De$ (the locus of the
critical values)
is a dianalytic curve with nodes and cusps as only singularities.
In this way  we may however have introduced also negative nodes, i.e.,
nodes which in local holomorphic coordinates are defined by the equation
$$ ( y -\bar x) ( y + \bar x) = 0. $$

Now, the projection onto the first factor $ \PP^1$ gives two 
monodromy factorizations:

\begin{df}
The {\bf vertical braid monodromy factorization} of a bidouble cover is
the factorization of the identity in the braid
group $\br_n $ of the sphere $ \PP^1$ associated to the first projection onto
  $ \PP^1$ of the pair $(Q, \De)$.
Here $n = 4 (b+d)$ is the vertical degree of $\De$, and the monodromy 
factorization
describes the movement of $n = 4 (b+d)$
points in a fibre $ F_0 = \PP^1$.
\end{df}

The Braid Monodromy factorization
  has factors in
the braid group of the sphere $\br_n$:
$$
\br_n \quad = \quad
\left\langle \: \s_1,..., \s_{n-1}\:\left|
\begin{array}{l}
\s_i\s_j = \s_j\s_i,\quad \text{ if } |i-j|>1\\
\s_i\s_{i+1}\s_i = \s_{i+1}\s_i\s_{i+1}\\
\s_1\cdots\s_{n-1}\s_{n-1}\cdots\s_1=1
\end{array}
\right.\right\rangle
$$

The factors are of three following types

\begin{enumerate}
\item
each cusp of $\De$ (which has non vertical tangent) yields a factor 
which is a conjugate of $\s_1^3$
\item
each vertical tangency (this takes place in smooth points only) 
yields a factor which is a conjugate of $\s_1$
\item
each trivial singularity yields a factor which is a conjugate of 
$\s_1^2$ in the case where we have a   positive node,
and  a factor which is a conjugate of $\s_1^{-2}$ in the case where 
we have a negative node.

\end{enumerate}

We already know, by Kas' theorem,  that the differentiable type
of the  Lefschetz fibration $F:  S \ra \PP^1$ given by  the
composition of the perturbed covering $\pi : S \ra Q$
with the  first projection onto $ \PP^1$ is encoded in the
M-equivalence class of the factorization of the identity in
the Mapping class group  $\MMM ap_g$  corresponding to the monodromy
factorization of this Lefschetz fibration (here $g$ is the genus of 
the base fibre $C_0$ of $F$) .

The second factorization is the homomorphic image of the
Braid Monodromy factorization.

It is obtained since the
braid monodromy factors are contained in a subgroup of
the Braid group $\br_n$, namely, the subgroup of liftable
braids, corresponding to diffeomorphisms which lift to the
fibre $C_0 $ of the Lefschetz fibration $F:  S \ra \PP^1$
which lies over $F_0$.

More precisely, we have the  covering $C_0 \setminus F^{-1} (\De) \ra 
F_0  \setminus \De$,
and  let $\mu : \pi_1 ( F_0  \setminus \De) \ra \gS_4$ be its 
monodromy homomorphism:
then the braid monodromy factors are contained in the subgroup
$$ \sL : = \{ \sigma \in \br_n | \mu \sigma_* = \mu \}  .$$

And each braid in $\sL$ determines a diffeomorphism of the pair  $( 
C_0 , C_0 \cap \De ) $
which in turn gives us an element in  $\MMM ap_g$ .

The easy but important observation is that the homomorphism
$$\sL \ra \MMM ap_g $$
is far from being injective.

In fact, the factors of the Braid Monodromy factorization  map to:

\begin{enumerate}
\item
each factor corresponding to a cusp of $\De$ maps to the identity
\item
each factor coming from a vertical tangency maps to a Dehn twist 
which is the lift of a conjugate  of $\s_1$
\item
each factor coming from a trivial singularity , be it  a   positive 
or a negative node,
maps to the identity.

\end{enumerate}

Before we state our characterization of the Braid Monodromy factorization,
it will be convenient to briefly describe the geometric picture of 
the degree $4$
covering  $F : C_0 \ra F_0$.

In the case where $b=d$ we may  start with an extremely   symmetric 
picture, taking the  bidouble cover $F'_0$ of
$\PP^1_{\C}$ of equation
\begin{itemize}
\item
$ z ^2 = f(y) , w  ^2 = g(y) : =  f (-y)$\\
and of the automorphism $\psi$  given by
\item
$ \psi (y,z,w) : = ( -y , w , z) $
\end{itemize}

and we see immediately that the curve $C_0$ has an automorphism group
generated by the previously mentioned automorphism $\psi$ and by the 
Galois group of the covering,
generated by
$$  g_1  (y,z,w) : = ( y , -z , w) , g_2  (y,z,w) : = ( y , z , - w).$$

This group is indeed generated by $g_1$ and by  the order $4$ element $ g_4 (y,z,w) 
: = ( -y , w , -z) $;
 since $g_1 g_4  g_1 g_4 = 1$, it is the dihedral group $D_4$.

We recall here once more the role of the automorphism $\psi$,
  namely that the following was proven in \cite{c-w}  :

\begin{prop}\label{mon}
The monodromy of $S$ over the unit circle $ \{ x | |x| = 1\}$ is trivial,
and the pair $(C_0, \psi)$ yields  the fibre over $ x=1$, considered as
     a differentiable  $2$-manifold, together with the isotopy class
of the  map yielding the twisted fibre sum corresponding to $S'$ .
\end{prop}

In the general case, perturbing the equations of $C_0$ amounts to 
making the degree 4 morphism
generic. In this context we can view the effect of the perturbation 
as replacing
the branch locus
$$ B_0':=  \{ \nodalf(y) = 0\}\cup \{ \nodalg(y) = 0\} = \{ 1, \dots 
,2b\}  \cup \{ - 1, \dots ,
-2d\}$$ with a new `doubled' branch locus
  $$ B_0:=  \{ f(y)g(y) = 0 \} = \{ 1\pm \epsilon \sqrt {-1}, \dots ,2b \pm
\epsilon \sqrt {-1}\} \cup \{ - 1 \pm \epsilon \sqrt {-1}, \dots , 
-2d \pm \epsilon \sqrt {-1}\}.$$

\newcommand{\quadrat}{\makebox(0,0){${\scriptscriptstyle \square}$}}
\newcommand{\kreis}{\makebox(0,0){$\bullet$}}
\unitlength=.53mm
\linethickness{.8pt}

\begin{figure}[htb]
\centering
\begin{picture}(240,240)(-140,-160)

\qbezier(-140,0)(-140,20)(-20,20)
\qbezier(100,0)(100,20)(20,20)
\qbezier(-140,0)(-140,-20)(-35,-20)
\qbezier(100,0)(100,-20)(35,-20)

\qbezier(-110,-5)(-105,-5)(-100,2)
\qbezier(-110,-5)(-115,-5)(-120,2)
\qbezier(-117,0)(-110,5)(-103,0)

\qbezier(-70,-5)(-65,-5)(-60,2)
\qbezier(-70,-5)(-75,-5)(-80,2)
\qbezier(-77,0)(-70,5)(-63,0)

\qbezier(70,-5)(65,-5)(60,2)
\qbezier(70,-5)(75,-5)(80,2)
\qbezier(77,0)(70,5)(63,0)

\qbezier(-3,3.5)(0,4)(3,3.5)
\qbezier(-38,0)(-30,3)(-18,3.5)
\qbezier(18,3.5)(30,3)(38,0)

\qbezier(-140,60)(-140,80)(0,80)
\qbezier(100,60)(100,80)(0,80)
\qbezier(-140,60)(-140,40)(-40,40)
\qbezier(100,60)(100,40)(40,40)

\qbezier(-110,55)(-105,55)(-100,62)
\qbezier(-110,55)(-115,55)(-120,62)
\qbezier(-117,60)(-110,65)(-103,60)

\qbezier(-70,55)(-65,55)(-60,62)
\qbezier(-70,55)(-75,55)(-80,62)
\qbezier(-77,60)(-70,65)(-63,60)

\qbezier(70,55)(65,55)(60,62)
\qbezier(70,55)(75,55)(80,62)
\qbezier(77,60)(70,65)(63,60)

\qbezier(-38,60)(0,65)(38,60)

\qbezier(40,40)(20,40)(0,8)
\qbezier(-35,-20)(-17.5,-20)(0,8)

\qbezier(-45,2)(-25,-8)(-20,0)
\qbezier(-20,0)(10,50)(45,62)

\put(100,-.5){\quadrat}
\put(77.5,-.5){\quadrat}
\put(62.5,-.5){\quadrat}
\put(39,-.5){\quadrat}
\put(100,59.5){\quadrat}
\put(77.5,59.5){\quadrat}
\put(62.5,59.5){\quadrat}
\put(39,59.5){\quadrat}

\put(-102.5,-.5){\kreis}
\put(-117.5,-.5){\kreis}
\put(-140,-.5){\kreis}
\put(-62.5,-.5){\kreis}
\put(-77.5,-.5){\kreis}
\put(-39,-.5){\kreis}
\put(-102.5,59.5){\kreis}
\put(-117.5,59.5){\kreis}
\put(-140,59.5){\kreis}
\put(-62.5,59.5){\kreis}
\put(-77.5,59.5){\kreis}
\put(-39,59.5){\kreis}

\put(105,-8){\makebox(0,0){\scriptsize \sf B$'_4$}}
\put(82,-8){\makebox(0,0){\scriptsize \sf B$'_3$}}
\put(61,-8){\makebox(0,0){\scriptsize \sf B$'_2$}}
\put(43,-8){\makebox(0,0){\scriptsize \sf B$'_1$}}
\put(105,67){\makebox(0,0){\scriptsize \sf B$''_4$}}
\put(82,67){\makebox(0,0){\scriptsize \sf B$''_3$}}
\put(61,67){\makebox(0,0){\scriptsize \sf B$''_2$}}
\put(43,67){\makebox(0,0){\scriptsize \sf B$''_1$}}

\put(-41,-8){\makebox(0,0){\scriptsize \sf D$'_1$}}
\put(-59,-8){\makebox(0,0){\scriptsize \sf D$'_2$}}
\put(-79,-8){\makebox(0,0){\scriptsize \sf D$'_3$}}
\put(-99,-8){\makebox(0,0){\scriptsize \sf D$'_4$}}
\put(-119,-8){\makebox(0,0){\scriptsize \sf D$'_5$}}
\put(-143,-8){\makebox(0,0){\scriptsize \sf D$'_6$}}
\put(-41,67){\makebox(0,0){\scriptsize \sf D$''_1$}}
\put(-59,67){\makebox(0,0){\scriptsize \sf D$''_2$}}
\put(-79,67){\makebox(0,0){\scriptsize \sf D$''_3$}}
\put(-99,67){\makebox(0,0){\scriptsize \sf D$''_4$}}
\put(-119,67){\makebox(0,0){\scriptsize \sf D$''_5$}}
\put(-143,67){\makebox(0,0){\scriptsize \sf D$''_6$}}

\put(8,19){\color{white}\circle*{5}}
\put(-8,19){\color{white}\circle*{5}}
\put(0,8){\color{white}\circle*{5}}
\put(0,29){\color{white}\circle*{5}}
\put(4,13.5){\color{white}\circle*{3}}
\put(-4,24){\color{white}\circle*{3}}

\qbezier(-40,40)(-20,40)(0,8)
\qbezier(35,-20)(17.5,-20)(0,8)

\qbezier(45,2)(25,-8)(20,0)
\qbezier(20,0)(-10,50)(-45,62)


\qbezier(-140,-70)(-140,-50)(-20,-50)
\qbezier(100,-70)(100,-50)(20,-50)
\qbezier(-140,-70)(-140,-90)(-35,-90)
\qbezier(100,-70)(100,-90)(35,-90)
\qbezier(-35,-90)(0,-90)(35,-90)
\qbezier(-20,-50)(0,-50)(20,-50)

\put(100,-70){\quadrat}
\put(77.5,-70){\quadrat}
\put(62.5,-70){\quadrat}
\put(39,-70){\quadrat}

\put(-102.5,-70){\kreis}
\put(-117.5,-70){\kreis}
\put(-140,-70){\kreis}
\put(-62.5,-70){\kreis}
\put(-77.5,-70){\kreis}
\put(-39,-70){\kreis}

\put(105,-78){\makebox(0,0){\scriptsize \sf B$_4$}}
\put(79,-78){\makebox(0,0){\scriptsize \sf B$_3$}}
\put(63.5,-78){\makebox(0,0){\scriptsize \sf B$_2$}}
\put(40,-78){\makebox(0,0){\scriptsize \sf B$_1$}}

\put(-38,-78){\makebox(0,0){\scriptsize \sf D$_1$}}
\put(-61,-78){\makebox(0,0){\scriptsize \sf D$_2$}}
\put(-76,-78){\makebox(0,0){\scriptsize \sf D$_3$}}
\put(-101,-78){\makebox(0,0){\scriptsize \sf D$_4$}}
\put(-116,-78){\makebox(0,0){\scriptsize \sf D$_5$}}
\put(-143,-78){\makebox(0,0){\scriptsize \sf D$_6$}}

{
\linethickness{1.2pt}
\qbezier(-80,-26)(-80,-36)(-80,-46)
\qbezier(-80,-46)(-78,-43)(-76,-40)
\qbezier(-80,-46)(-82,-43)(-84,-40)
}
\put(-105,-38){\makebox(0,0){Galois}}


\qbezier(-140,-140)(-140,-120)(-20,-120)
\qbezier(100,-140)(100,-120)(20,-120)
\qbezier(-140,-140)(-140,-160)(-35,-160)
\qbezier(100,-140)(100,-160)(35,-160)
\qbezier(-35,-160)(0,-160)(35,-160)
\qbezier(-20,-120)(0,-120)(20,-120)

\put(99.5,-143){\quadrat}
\put(77.5,-143){\quadrat}
\put(62.5,-143){\quadrat}
\put(39,-143){\quadrat}
\put(99.5,-137){\quadrat}
\put(77.5,-137){\quadrat}
\put(62.5,-137){\quadrat}
\put(39,-137){\quadrat}

\put(-102.5,-143){\kreis}
\put(-117.5,-143){\kreis}
\put(-139.5,-143){\kreis}
\put(-62.5,-143){\kreis}
\put(-77.5,-143){\kreis}
\put(-39,-143){\kreis}
\put(-102.5,-137){\kreis}
\put(-117.5,-137){\kreis}
\put(-139.5,-137){\kreis}
\put(-62.5,-137){\kreis}
\put(-77.5,-137){\kreis}
\put(-39,-137){\kreis}

\put(105,-148){\makebox(0,0){\scriptsize \sf B$'_4$}}
\put(79,-149){\makebox(0,0){\scriptsize \sf B$'_3$}}
\put(63.5,-149){\makebox(0,0){\scriptsize \sf B$'_2$}}
\put(40,-149){\makebox(0,0){\scriptsize \sf B$'_1$}}

\put(-38,-149){\makebox(0,0){\scriptsize \sf D$'_1$}}
\put(-61,-149){\makebox(0,0){\scriptsize \sf D$'_2$}}
\put(-76,-149){\makebox(0,0){\scriptsize \sf D$'_3$}}
\put(-101,-149){\makebox(0,0){\scriptsize \sf D$'_4$}}
\put(-116,-149){\makebox(0,0){\scriptsize \sf D$'_5$}}
\put(-143,-148){\makebox(0,0){\scriptsize \sf D$'_6$}}

\put(105,-131){\makebox(0,0){\scriptsize \sf B$''_4$}}
\put(79,-131){\makebox(0,0){\scriptsize \sf B$''_3$}}
\put(63.5,-131){\makebox(0,0){\scriptsize \sf B$''_2$}}
\put(40,-131){\makebox(0,0){\scriptsize \sf B$''_1$}}

\put(-38,-131){\makebox(0,0){\scriptsize \sf D$''_1$}}
\put(-61,-131){\makebox(0,0){\scriptsize \sf D$''_2$}}
\put(-76,-131){\makebox(0,0){\scriptsize \sf D$''_3$}}
\put(-101,-131){\makebox(0,0){\scriptsize \sf D$''_4$}}
\put(-116,-131){\makebox(0,0){\scriptsize \sf D$''_5$}}
\put(-143,-131){\makebox(0,0){\scriptsize \sf D$''_6$}}

{
\linethickness{1.2pt}
\qbezier(20,-26)(20,-36)(20,-46)
\qbezier(20,-96)(20,-106)(20,-116)
\qbezier(20,-116)(18,-113)(16,-110)
\qbezier(20,-116)(22,-113)(24,-110)
}
\put(68,-105){\makebox(0,0){perturbed}}

\end{picture}

\end{figure}

Before giving the full vertical Braid monodromy factorization, we 
show a simpler byproduct of our
description of the Braid monodromy: it  makes effective and precise 
the  contents of the above intuitive picture.

The theorem is stated in its most understandable form, namely a 
pictorial one, recalling  that
any segment in the plane determines a (non standard) half twist in 
$\br_n, \ n = 4 (b+d)$.

\begin{thm}
\label{bmf}
There is a braid monodromy factorization of the curve $\De$
associated to a bidouble cover $S$ of type $a,b,c,d$
       whose braid monodromy group $H\subset\br_{4(b+d)}$ is generated,
unless we are in the cases  $$ {\rm(I) }\  c = 2a \ {\rm and } \  d =
2b \ , \ {\rm or\  (II)} \ a = 2c
\ {\rm and } \  b = 2d \,$$

     by
the following powers of half-twists:
$\s_{a_i},\s_{c_i}$ for $i=1,\dots 2b-1$, $\s_{b_\imath},\s_{d_\imath}$,
for $\imath=1,\dots,2d-1$,
$\s_{p_{2b}}^2,\s_{q_{2d}}^2,\s_{s},\s_{u'}^3,\s_{u''}^3$.

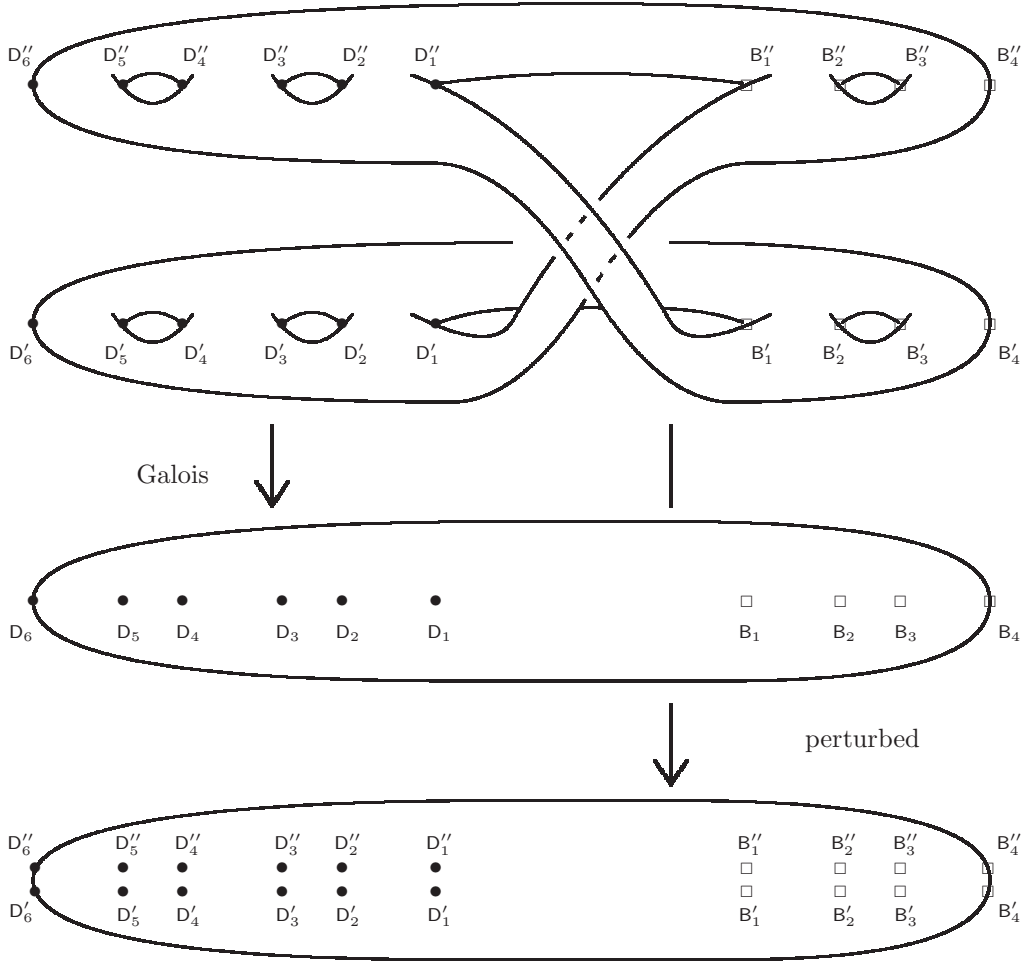
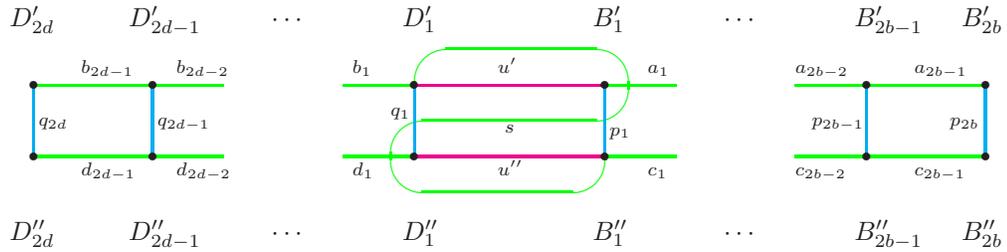
\unitlength=.9pt
\begin{picture}(420,120)(-210,-60)

{\color{green}
\put(-200,-15){\line(1,0){80}}
\put(-70,-15){\line(1,0){30}}
\put(200,-15){\line(-1,0){80}}
\put(70,-15){\line(-1,0){30}}
\put(-200,15){\line(1,0){80}}
\put(-70,15){\line(1,0){30}}
\put(200,15){\line(-1,0){80}}
\put(70,15){\line(-1,0){30}}

\put(5,15){\oval(90,30)[t]}
\put(0,15){\oval(100,30)[rb]}
\put(0,-15){\oval(100,30)[lt]}
\put(-5,-15){\oval(90,30)[b]}
\color{magenta}
\put(-40,15){\line(1,0){80}}
\put(-40,-15){\line(1,0){80}}
\color{cyan}
\put(-200,-15){\line(0,1){30}}
\put(-150,-15){\line(0,1){30}}
\put(-40,-15){\line(0,1){30}}
\put(40,-15){\line(0,1){30}}
\put(150,-15){\line(0,1){30}}
\put(200,-15){\line(0,1){30}}
}

\put(40,15){\circle*{3}}
\put(150,15){\circle*{3}}
\put(200,15){\circle*{3}}
\put(-40,15){\circle*{3}}
\put(-150,15){\circle*{3}}
\put(-200,15){\circle*{3}}
\put(40,-15){\circle*{3}}
\put(150,-15){\circle*{3}}
\put(200,-15){\circle*{3}}
\put(-40,-15){\circle*{3}}
\put(-150,-15){\circle*{3}}
\put(-200,-15){\circle*{3}}

\put(-210,-50){$D''_{2d}$}
\put(-160,-50){$D''_{2d-1}$}
\put(-45,-50){$D''_{1}$}
\put(35,-50){$B''_{1}$}
\put(145,-50){$B''_{2b-1}$}
\put(190,-50){$B''_{2b}$}
\put(-210,40){$D'_{2d}$}
\put(-160,40){$D'_{2d-1}$}
\put(-45,40){$D'_{1}$}
\put(35,40){$B'_{1}$}
\put(145,40){$B'_{2b-1}$}
\put(190,40){$B'_{2b}$}

\put(-100,-50){$\cdots$}
\put(-100,40){$\cdots$}
\put(90,-50){$\cdots$}
\put(90,40){$\cdots$}

\put(-198,-2){$\scriptstyle q_{2d}$}
\put(-148,-2){$\scriptstyle q_{2d-1}$}
\put(-50,2){$\scriptstyle q_{1}$}
\put(42,-6){$\scriptstyle p_{1}$}
\put(127,-2){$\scriptstyle p_{2b-1}$}
\put(185,-2){$\scriptstyle p_{2b}$}

\put(-180,20){$\scriptstyle b_{2d-1}$}
\put(-140,20){$\scriptstyle b_{2d-2}$}
\put(-66,20){$\scriptstyle b_{1}$}
\put(-5,20){$\scriptstyle u'$}
\put(58,20){$\scriptstyle a_{1}$}
\put(120,20){$\scriptstyle a_{2b-2}$}
\put(170,20){$\scriptstyle a_{2b-1}$}
\put(-180,-23){$\scriptstyle d_{2d-1}$}
\put(-140,-23){$\scriptstyle d_{2d-2}$}
\put(-66,-23){$\scriptstyle d_{1}$}
\put(-5,-23){$\scriptstyle u''$}
\put(58,-23){$\scriptstyle c_{1}$}
\put(120,-23){$\scriptstyle c_{2b-2}$}
\put(170,-23){$\scriptstyle c_{2b-1}$}

\put(-2,-5){$\scriptstyle s$}
\end{picture}\\
The factorization is such that
\begin{enumerate}
\item
each ($\pm$)full-twist factor is of type $p$ or $q$,
\item
the weighted count of ($\pm$) full-twist factors of type $p$ yields
$8ab-2(ad+bc)$,
\item
the weighted count of ($\pm$) full-twist factors of type $q$ yields
$8cd-2(ad+bc)$.
\end{enumerate}
\end{thm}

The result  relies on a complete description of the braid monodromy
factorization class
associated to $\De$: this is given in Theorem \ref{BMF}, which proves
indeed much more than what we need
for the present purposes.

The most important application  is  the following Main Theorem:

\begin {teo}

Assume that $S$ is an $abcd$-surface, $S'$ is an $a'b'c'd'$-surface,
that $S$ and $S'$ are homeomorphic, and that their  vertical Braid 
monodromy factorizations
are stably equivalent. Then, up to swapping  the pair $(a,b)$ with the pair $(c,d)$, we have
$$ a+c = a' + c' , b+d = b' + d'  , ab = a'b' , cd = c' d' .$$
The vertical Braid monodromy factorizations associated to
an abc-surface $S$  and to
         an $a'b'c'$-surface  $S'$ are not stably equivalent, except in the
trivial cases
         $ a = a', b = b', c =  c'$ or $ a = c', b = b', c =  a'$.
         \end{teo}

\proof
We sketch only the first calculation, which is not contained in \cite{clw}.

Since $S$ and $S'$ are homeomorphic, they have the same $\chi$ and $K^2$,
thus first of all $ ab + cd = a'b' + c'  d'$. By the previous theorem 
\ref{bmf} and the forthcoming
proposition (\ref{nonconj}) that the factors of type $p$ are not 
stably conjugate to the factors of
type
$q$, it follows that
$ab = a'b' , cd = c' d' .$

Then by the same formulae of  theorem \ref{bmf}  $ (a+c)(b+d) = 
(a'+c')(b'+d')$,
hence, since $K^2_S = (a+c-2) (b+d-2)$, we get $ (a+c) + (b+d) = 
(a'+c') + (b'+d')$,
and therefore $ a+c = a' + c' , b+d = b' + d' $, up to swapping the pairs.

\QED

We start now to prepare the notation for the factorization theorem.

\begin{prop}
\label{c-mono}
Up to homotopy in small discs containing pairs $D'_i,D''_i$ or $B_\jmath',
B_\jmath''$ the punctured fibre $F_0$ is given by the following picture
\\
\unitlength=.9pt
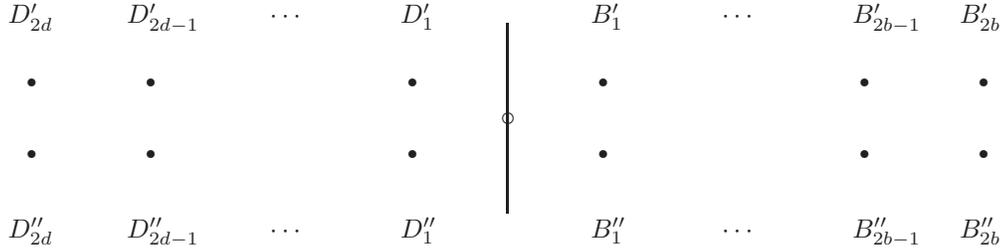
\begin{picture}(420,120)(-210,-60)

\put(0,0){\circle{5}}
\put(0,-40){\line(0,1){80}}

\put(40,15){\circle*{3}}
\put(150,15){\circle*{3}}
\put(200,15){\circle*{3}}
\put(-40,15){\circle*{3}}
\put(-150,15){\circle*{3}}
\put(-200,15){\circle*{3}}
\put(40,-15){\circle*{3}}
\put(150,-15){\circle*{3}}
\put(200,-15){\circle*{3}}
\put(-40,-15){\circle*{3}}
\put(-150,-15){\circle*{3}}
\put(-200,-15){\circle*{3}}

\put(-210,-50){$D''_{2d}$}
\put(-160,-50){$D''_{2d-1}$}
\put(-45,-50){$D''_{1}$}
\put(35,-50){$B''_{1}$}
\put(145,-50){$B''_{2b-1}$}
\put(190,-50){$B''_{2b}$}
\put(-210,40){$D'_{2d}$}
\put(-160,40){$D'_{2d-1}$}
\put(-45,40){$D'_{1}$}
\put(35,40){$B'_{1}$}
\put(145,40){$B'_{2b-1}$}
\put(190,40){$B'_{2b}$}

\put(-100,-50){$\cdots$}
\put(-100,40){$\cdots$}
\put(90,-50){$\cdots$}
\put(90,40){$\cdots$}

\end{picture}\\
and the covering monodromy $\mu$ of the perturbed bidouble cover,
with respect to the origin,
is given by
$$
\begin{array}{rcl@{\hspace{2cm}}rcl}
\oo_{D'_j} & \mapsto & (12), & \oo_{B'_j} & \mapsto & (13),\\
\oo_{D''_j} & \mapsto & (34), &
\oo_{B''_j} & \mapsto & (24),\\
\end{array}
$$
where $\oo_P$ is any simple closed path around $P$
not crossing the imaginary axis $i\RR\subset F_0$.
\end{prop}

\proof
As we remarked before, in the unperturbed Galois cover case we have
real points $D_i$,
$B_\jmath$, only, the  $D_i$'s with positive real coordinate,
the  $B_j$'s with negative real coordinates.
And the covering's  monodromy is $(12)(34)$, resp.\ $(13)(24)$
for paths which do not cross the imaginary axis.

The deformation then splits each branch point into two, hence the
corresponding monodromies (which are transpositions)
must be $(12)$ and $(34)$, resp.\ $(13)$
and
$(24)$ (which commute).
After a suitable homotopy the points are in the positions  given by
the picture.

\qed

Referring to  the  above figures and to the figure below we introduce the following notation for
a system of arcs, which
are uniquely determined (up to homotopy) by their endpoints and the
property that
they are monotonous in the real coordinate and do not pass below
any puncture; i.e.,  if they share the real coordinate with a puncture
they have larger imaginary coordinate.
$$
\begin{array}{ll@{\hspace*{2cm}}ll}
p_{i}: & B'_i, B''_i, &
q_{\imath}: & D'_\imath, D''_\imath,\\
a_{ij}: & B'_i, B'_j, &
b_{\imath\jmath}: & D'_ \imath, D'_\jmath, \\
c_{ij}: & B''_i, B''_j, &
d_{\imath\jmath}: & D''_\imath, D''_\jmath, \\[2mm]
u'_{i\jmath}: & B'_i, D'_\jmath &
u''_{i\jmath}: & B''_i, D''_\jmath
\end{array}
$$

Further, there are arcs $s_{ij}$ connecting $D'_i$ with $B''_j$
according to a pattern which is shown in the next  figure \ref{arcs3}.

\unitlength=.92pt
\begin{figure} 
\begin{picture}(420,120)(-205,-60)

\put(40,5){\circle*{3}}
\put(135,5){\circle*{3}}
\put(200,5){\circle*{3}}
\put(-40,5){\circle*{3}}
\put(-135,5){\circle*{3}}
\put(-200,5){\circle*{3}}
\put(40,-20){\circle*{3}}
\put(135,-20){\circle*{3}}
\put(200,-20){\circle*{3}}
\put(-40,-20){\circle*{3}}

\put(-135,-20){\circle*{3}}
\put(-200,-20){\circle*{3}}

\put(-75,-20){\circle*{3}}
\put(-75,5){\circle*{3}}

\put(75,-20){\circle*{3}}
\put(75,5){\circle*{3}}

\put(-205,-50){$D''_{2d}$}
\put(-140,-50){$D''_{i}$}
\put(-80,-50){$D''_{2}$}
\put(-45,-50){$D''_{1}$}
\put(35,-50){$B''_{1}$}
\put(70,-50){$B''_{2}$}
\put(130,-50){$B''_{j}$}
\put(190,-50){$B''_{2b}$}
\put(-205,40){$D'_{2d}$}
\put(-140,40){$D'_{i}$}
\put(-45,40){$D'_{1}$}
\put(-80,40){$D'_{2}$}
\put(35,40){$B'_{1}$}
\put(70,40){$B'_{2}$}
\put(130,40){$B'_{j}$}
\put(190,40){$B'_{2b}$}

\put(-110,-50){$\cdots$}
\put(-110,40){$\cdots$}
\put(100,-50){$\cdots$}
\put(100,40){$\cdots$}

\put(-170,-50){$\cdots$}
\put(-170,40){$\cdots$}
\put(160,-50){$\cdots$}
\put(160,40){$\cdots$}

\qbezier(-135,5)(-122,21)(-90,26)
\qbezier(-90,26)(-30,35)(40,13)
\qbezier(40,13)(47,10,)(47,5)
\qbezier(47,5)(47,-3)(41,-3)
\qbezier(41,-3)(37,-3,)(32,0)
\qbezier(32,0)(-10,23)(-60,23)
\qbezier(-60,23)(-95,23)(-120,10)
\qbezier(-120,10)(-140,-2)(-141,-20)
\qbezier(-141,-20)(-141.5,-27)(-135,-27)
\qbezier(-135,-27)(-130,-27)(-128,-18)
\qbezier(-128,-18)(-124,-3)(-110,5)
\qbezier(40,-20)(-40,40)(-110,5)

\qbezier(135,-20)(120,18)(50,18)
\qbezier(50,18)(-20,18)(-40,-20)

\qbezier(135,5)(100,25)(50,25)
\qbezier(50,25)(0,25)(-40,5)

\put(105,-5){$\scriptstyle u''_{1j}$}
\put(115,20){$\scriptstyle u'_{1j}$}
\put(-110,-5){$\scriptstyle s_{i1}$}
\end{picture}

\caption{arcs for the theorem}
 \label{arcs3}
\end{figure}
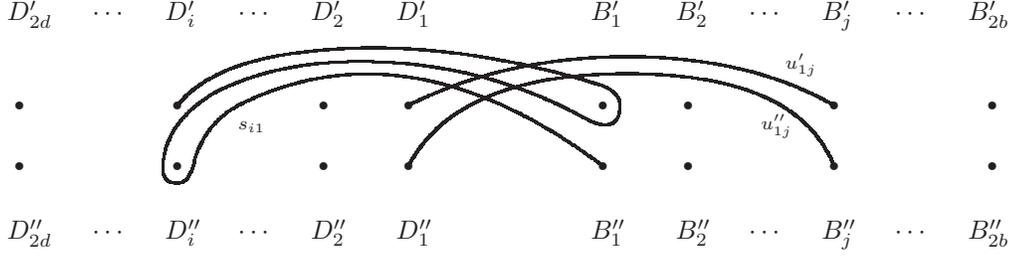

\begin{thm}
\label{BMF}
The factorization is a product of factorizations
$$
\begin{matrix}
\underbrace{
\bigg(\b_f \circ
\underbrace{(\s_{p_1}^{\pm2} \circ \cdots \circ
\s_{p_1}^{\pm2})}_{|2b-d|\text{ factors}}
\circ \b_{fg} \bigg)\circ \cdots}_{2a\text{ repetitions}}
\circ
\underbrace{(\s_{p_1}^{\pm2} \circ \cdots \circ
\s_{p_{2b}}^{\pm2})\circ \cdots}_{|2a-c|\text{ repetitions}}
\\[2mm]
\hspace*{30mm}\circ
\underbrace{(\s_{q_1}^{\pm2} \circ \cdots \circ
\s_{q_{2d}}^{\pm2})\circ \cdots}_{|2c-a|\text{ repetitions}}
\circ
\underbrace{
\bigg(\b_g \circ
\underbrace{(\s_{q_1}^{\pm2} \circ \cdots \circ
\s_{q_1}^{\pm2})}_{|2d-b|\text{ factors}}
\circ \b_{gf} \bigg)\circ \cdots}_{2c\text{ repetitions}}
\end{matrix}
$$
where the sign of the exponents $2$ is constant inside a pair of
brackets and is the sign of the number which determines the number
of factors, ie.\ $(2b-d)$ resp.\ $(2a-c)$, $(2c-a)$ or $(2d-b)$,
and
where the $\b$'s further decompose as products of factorizations
\[\begin{matrix}
\b_f & = & \b_{f,2}\circ \cdots \circ \b_{f,2b},\quad
\b_{fg} & = & \b_{fg,2d}\circ \cdots \circ \b_{fg,1},\\
\b_{g} & = & \b_{g,2} \circ \cdots \circ \b_{g,2d},\quad
\b_{gf} & = & \b_{gf,2b} \circ \cdots \circ \b_{gf,1},
        \end{matrix}\]
based on elementary factorizations each having four factors
\[\begin{matrix}
\b_{f,i} & = & \s_{a_{1i}}\circ \s_{p_1}^2\s_{c_{1i}}\s_{p_1}^{-2}\circ
\s_{a_{1i}}\circ \s_{p_1}^2\s_{c_{1i}}\s_{p_1}^{-2},\qquad
\b_{fg,j} & = & \s_{u_{1j}}^3\circ \s_{s_{1j}}\circ
\s_{u_{1j}'}^3\circ \s_{u_{1j}''}^3,\\
\b_{g,j} & = & \s_{b_{1j}}\circ \s_{q_1}^2\s_{d_{1j}}\s_{q_1}^{-2}\circ
\s_{b_{1j}}\circ \s_{q_1}^2\s_{d_{1j}}\s_{q_1}^{-2},\qquad
\b_{gf,i} & = & \s_{u_{i1}}^3\circ \s_{s_{i1}}\circ
\s_{u_{i1}'}^3\circ \s_{u_{i1}''}^3.
        \end{matrix}\]
The elementary factorizations originate in the regeneration of nodes
of the branch curve of the corresponding bidouble Galois-cover.
\end{thm}

While it is complicated to explain in detail the above theorem,
we would like to point out some other  important  tools in the local 
analysis of
`regenerations' of the singularities of the branch curve of the bidouble cover.

In the  case of a `proper node'  ( i.e., a point where $f=g=0$) we 
get a cusp-cluster, four critical
points of which
three are cusps and the last is a vertical tangency point.

The braid monodromy for the regeneration into a cusp-cluster
had been thoroughly investigated before in \cite{c-w2},
where the braid monodromy factorization type  had been
determined up to
Hurwitz equivalence and simultaneous conjugation.

In \cite{clw} we were able to show how one can determine the
      Hurwitz equivalence class of the factorization
from the datum of  the product of the
factors.

\newcommand{\braid}{\cg}

\begin{prop}
\label{hur1}
The braid monodromy factorization of a (regenerated) cusp-cluster with
product $\s_2^3\s_1\s_3\s_2\s_1^2\s_3^2$ is
Hurwitz-equivalent to the factorization
$$
\s_2^3 \circ  \s_1 \s_3 \s_2 \s_3\inv \s_1\inv \circ \s_1^3 \circ \s_3^3.
$$
\end{prop}

In the case of an improper node (i.e., a node of $\nodalf = 0$ or of 
$\nodalg = 0$) we get a cluster of vertical tangents, four
critical points, all of which are vertical tangency points.

\begin{prop}
\label{vt-fact}
The braid monodromy factorization associated to a cluster of
tangents corresponding to an improper node  is given by
\[
\s_2\s_3\s_2\inv \circ \s_1\s_2\s_1\inv \circ
\s_2\s_3\s_2\inv \circ \s_1\s_2\s_1\inv
\]
up to Hurwitz equivalence and simultaneous conjugation.
\end{prop}

\section{Another proof
of the non conjugacy theorem}

\unitlength=.9pt
\begin{figure}
\begin{picture}(420,120)(-210,-60)

{\color{green}
\put(-200,-15){\line(1,0){80}}
\put(-70,-15){\line(1,0){30}}
\put(200,-15){\line(-1,0){80}}
\put(70,-15){\line(-1,0){30}}
\put(-200,15){\line(1,0){80}}
\put(-70,15){\line(1,0){30}}
\put(200,15){\line(-1,0){80}}
\put(70,15){\line(-1,0){30}}

\put(5,15){\oval(90,30)[t]}
\put(0,15){\oval(100,30)[rb]}
\put(0,-15){\oval(100,30)[lt]}
\put(-5,-15){\oval(90,30)[b]}
\color{magenta}
\put(-40,15){\line(1,0){80}}
\put(-40,-15){\line(1,0){80}}
\color{cyan}
\put(-200,-15){\line(0,1){30}}
\put(-150,-15){\line(0,1){30}}
\put(-40,-15){\line(0,1){30}}
\put(40,-15){\line(0,1){30}}
\put(150,-15){\line(0,1){30}}
\put(200,-15){\line(0,1){30}}
}

\put(40,15){\circle*{3}}
\put(150,15){\circle*{3}}
\put(200,15){\circle*{3}}
\put(-40,15){\circle*{3}}
\put(-150,15){\circle*{3}}
\put(-200,15){\circle*{3}}
\put(40,-15){\circle*{3}}
\put(150,-15){\circle*{3}}
\put(200,-15){\circle*{3}}
\put(-40,-15){\circle*{3}}

\put(-150,-15){\circle*{3}}
\put(-200,-15){\circle*{3}}

\put(-210,-50){$D''_{2d}$}
\put(-160,-50){$D''_{2d-1}$}
\put(-45,-50){$D''_{1}$}
\put(35,-50){$B''_{1}$}
\put(145,-50){$B''_{2b-1}$}
\put(190,-50){$B''_{2b}$}
\put(-210,40){$D'_{2d}$}
\put(-160,40){$D'_{2d-1}$}
\put(-45,40){$D'_{1}$}
\put(35,40){$B'_{1}$}
\put(145,40){$B'_{2b-1}$}
\put(190,40){$B'_{2b}$}

\put(-100,-50){$\cdots$}
\put(-100,40){$\cdots$}
\put(90,-50){$\cdots$}
\put(90,40){$\cdots$}

\put(-198,-2){$\scriptstyle
q_{2d}$}
\put(-148,-2){$\scriptstyle
q_{2d-1}$}
\put(-50,2){$\scriptstyle
q_{1}$}
\put(42,-6){$\scriptstyle p_{1}$}
\put(127,-2){$\scriptstyle
p_{2b-1}$}
\put(185,-2){$\scriptstyle
p_{2b}$}

\put(-180,20){$\scriptstyle
b_{2d-1}$}
\put(-140,20){$\scriptstyle
b_{2d-2}$}
\put(-66,20){$\scriptstyle
b_{1}$}
\put(-5,20){$\scriptstyle u_{x}$}
\put(58,20){$\scriptstyle
a_{1}$}
\put(120,20){$\scriptstyle
a_{2b-2}$}
\put(170,20){$\scriptstyle
a_{2b-1}$}
\put(-180,-23){$\scriptstyle
d_{2d-1}$}
\put(-140,-23){$\scriptstyle
d_{2d-2}$}
\put(-66,-23){$\scriptstyle
d_{1}$}
\put(-5,-23){$\scriptstyle u_{y}$}
\put(58,-23){$\scriptstyle
c_{1}$}
\put(120,-23){$\scriptstyle
c_{2b-2}$}
\put(170,-23){$\scriptstyle
c_{2b-1}$}

\put(-2,-5){$\scriptstyle
s$}
\end{picture}

\caption{isotoped generators}
\end{figure}
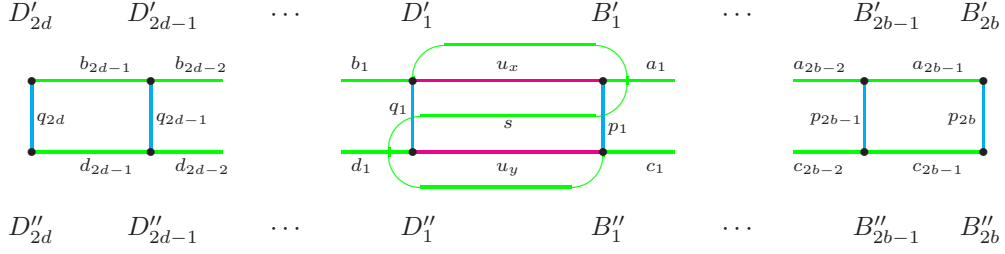

This section is based on the classical correspondence between
4-tuple covers and triple covers, given by the surjection
$\gS_4 \ra \gS_3$ whose kernel is the Klein group $(\Z / 2 )^2$.

This idea was also used in \cite{clw} to  obtain a triple cover $D_0 
\ra F_0$ corresponding to a quadruple cover $C_0 \ra F_0$ and a
resulting homomorphism of its  mapping class
group  to the symplectic group acting on the  $\Z / 2$  homology.

In this way in the article \cite{clw} we gave a proof that the 
elements $\sigma^2_p$
and $\sigma^2_q$ in the
braid monodromy group $H$ are not conjugate even in the stabilized
braid monodromy group $\hat H$.

Observe that in the braid group there is a unique square root $\sigma$ of
a full twist: since a full twist determines the homotopy class of an
arc between two punctures, and this arc determines
a half twist. In the braid group of the sphere, this square root
is not unique (as pointed out by a referee), since another square root
is also given by $\sigma \Delta^2$, where $ \Delta^2$ is a generator of the centre,
and has order $2$. However, there is a unique square root in the conjugacy class of half-twists.

    Hence it suffices
to show that there is no
$h$ in
$\hat H$ such that
$\sigma_p h= h \sigma_q$.

Our proof there exploited the action of the group $\hat H^+$
generated by $\hat H, \sigma_p$ and $\sigma_q$
on the $\mathbb Z/2\mathbb Z$ homology
of a the  triple cover $D_0$ of $ F_0  \cong \mathbb P^1$.

In the alternative proof we are going to present, we look instead at the
action of $\hat H^+$ on factorizations in $\gS_4$ of length $n=4b+4d$,
restriction of the Hurwitz action of the braid group $\br_n$.

In fact it suffices to consider the $\hat H^+$-orbit
$\mathbb O$ of the factorization
$$
\tau^0:\:\:
(12)\circ(34)\circ(12)\circ(34)\cdots(12)\circ(34)\circ(13)\circ(24)
\circ(13)\circ(24)\cdots(13)\circ(24),
$$
with $2b$ factors equal to $(12)$, respectively to  $(34)$, and
$2d$ factors  equal to $(13)$, respectively to $(24)$.
Note that $\tau^0$ is the factorization sequence corresponding to
the covering monodromy for the $4:1$ cover given by the
restriction of our bidouble cover to the vertical projective line
lying over the origin $x=0$.

The orbit $\mathbb O$ is rather small by the following lemma.

\begin{lemma}
The orbit $\mathbb O$ is contained in the set $\hat{\mathbb O}$ of
factorizations made of
transpositions
$$
    \tau_1\circ \tau_2\circ \dots\circ \tau_n
$$
such that for $i\leq 4b$ either $\tau_i=(12)$ or $\tau_i=(34)$, and
for $i>4b$ either $\tau_i=(13)$ or $\tau_i=(24)$.
\end{lemma}

\proof
We must show that all elements in $\hat H^+$ leave invariant the set
given in the claim.

To this purpose we exploit the surjection $\gS_4\to \gS_3$. It maps
transpositions to transpositions.
Naturally the induced map on factorizations commutes with the
Hurwitz action.

The upshot is now that the factorizations in the lemma are exactly
the factorizations whose factors are transpositions and which map
to one and the same factorization in $\gS_3$, namely
$$
\bar\tau^0:\:\:
(12)\circ(12)\circ(12)\circ\cdots\circ(12)\circ(13)\circ(13)
\circ(13)\circ\cdots\circ(13).
$$
    All generators of the braid monodromy subgroup $H$ act trivially
on the monodromy factorization in $\gS_4$,
hence also on the corresponding factorization $\bar\tau^0$  in $\gS_3$.
Likewise do the elements $\sigma_p,\sigma_q$: since for instance
$\sigma_q$ is the half twist connecting two neighbouring branch points
whose respective  local monodromies in $\gS_4$ yield consecutive
factors $(1,2)$ and
$(3,4)$: thus  $\sigma_q$ exchanges these two factors.

It remains to prove that also the additional elements needed
to generate $\hat H$ act trivially.

By definition of stable equivalence such additional generators are
full-twists $\sigma_\alpha^2$
which preserve  the monodromy factorization $\tau^0$ in $\gS_4$, hence
also  $\bar\tau^0$
and our claim follows.

\qed

We record now in detail how our generators of $\hat H^+$ do act.

\begin{lemma}
\label{triv-act}
The following braids act trivially on all elements in $\hat{\mathbb O}$:
$$
\sigma_{u'}^3,\sigma_{u''}^3,
$$
and all additional generators of $\hat H$ which are not in $H$.
\end{lemma}

\proof

Let us use as before the plane lexicographic ordering for the branch points.

Then one can easily verify that the element $\sigma_{4d}^3$ acts
trivially on all
factorizations, since $\tau_{4d}$ and $\tau_{4d+1}$ are always
non commuting transpositions.

Hence both elements $\sigma_{u'}^3,\sigma_{u''}^3$, which are  $H$-conjugate to
$\sigma_{4d}^3$,  act both trivially.

Consider now an additional generator $\sigma_\alpha^2$. It
preserves $\tau^0$ and may be written as
$\beta\sigma_1^2\beta^{-1}$.

Hence $\sigma_1^2$ acts trivially on $\bar\tau^0\beta$.
Therefore the first two factors of $\bar\tau^0\beta$ must commute.
Since both belong to $\gS_3$ they are even equal.
The corresponding factors in any lift of $\bar\tau^0\beta$
to a  factorization of transpositions in $\gS_4$ are then either equal
or disjoint, so $\sigma_1^2$ acts trivially also on all such lifts.

We conclude that $\sigma_\alpha^2=\beta\sigma_1^2\beta^{-1}$
acts trivially on all elements in $\hat{\mathbb O}$.

\qed

\begin{lemma}
\label{perm-act}
The following types of braids act by the corresponding transpositions
on  factorizations in $\hat{\mathbb O}$:
$$
\sigma_a, \sigma_c, \sigma_b, \sigma_d,
\sigma_p, \sigma_q.
$$
\end{lemma}

\proof
We observe that all elements $\sigma_i$ with $i\neq 4d$ act
by the corresponding transposition
on the  factorizations in $\hat{\mathbb O}$,
because consecutive transpositions commute except for the
transpositions in positions $4d$ and $4d+1$.

Hence the elements $\sigma_a, \sigma_c, \sigma_b, \sigma_d,
\sigma_p$ and $\sigma_q$,
which can be written as a product of factors different from
    $\sigma_{4d}$  act by the corresponding
transpositions.

\qed

\begin{lemma}
\label{s-act}
The action of the `snake' half twist $\sigma_s$ on  factorizations in
$\hat{\mathbb O}$
is given as follows:
\begin{enumerate}
\item
If the composition $\tau_{4d-1}\tau_{4d}\tau_{4d+1}\tau_{4d+2}$ is
trivial or equal to $\pi=(14)(23)$, then $\sigma_s$ acts trivially.
\item
Otherwise $\sigma_s$ acts by
$$
\tau_i \quad \mapsto \quad
\left\{
\begin{array}{cl}
\pi\tau_i\pi & \text{ if } 4d-1\leq i \leq 4d+2\\
\tau_i & \text{ else }
\end{array}
\right.
$$
\end{enumerate}
\end{lemma}

\proof
The `snake' half twist $\sigma_s$ is given in terms of the Artin generators as
$$
\sigma_s = \sigma_{4d}(\sigma_{4d-1}^{2}\sigma_{4d+1}^{2})\sigma_{4d} 
(\sigma_{4d+1}^{-2}\sigma_{4d-1}^{-2})\sigma_{4d}^{-1}
$$
Hence the element $\sigma_s$ acts on the subfactorization
$\tau_{4d-1}\circ\tau_{4d}\circ\tau_{4d+1}\circ\tau_{4d+2}$ and
leaves the remaining parts unchanged.
It thus suffices to know how $\sigma_s$ acts on all 16 possible
subfactorizations. 

First we record the action on four specific factorizations.
We decompose the action into five steps, each consisting of
one or two elementary Hurwitz moves.
In the first step we act by $\sigma_{4d}$, 
then by $(\sigma_{4d-1}^{2} \sigma_{4d+1}^{2})$, by $\sigma_{4d}$, 
by $(\sigma_{4d+1}^{-2}\sigma_{4d-1}^{-2})$ and finally by $\sigma_{4d}^{-1}$.
\begin{align*}
&\qquad (12)\circ(12)\circ(13)\circ(13) 
&\qquad (12)\circ(34)\circ(13)\circ(24)\\
& \mapsto\quad (12)\circ(23)\circ(12)\circ(13) 
&\mapsto\quad (12)\circ(14)\circ(34)\circ(24)\\
& \mapsto\quad (23)\circ(13)\circ(13)\circ(23) 
&\mapsto\quad (14)\circ(24)\circ(24)\circ(23)\\
& \mapsto\quad (23)\circ(13)\circ(13)\circ(23) 
&\mapsto\quad (14)\circ(24)\circ(24)\circ(23)\\
& \mapsto\quad (12)\circ(23)\circ(12)\circ(13) 
&\mapsto\quad(12)\circ(14)\circ(34)\circ(24)\\
& \mapsto\quad (12)\circ(12)\circ(13)\circ(13) 
&\mapsto\quad (12)\circ(34)\circ(13)\circ(24)
\end{align*}
We check that these first two factorizations belong to case (1) of the lemma
and remain indeed unchanged.
\begin{align*}
&\qquad (12)\circ(34)\circ(13)\circ(13) 
&\qquad (12)\circ(12)\circ(13)\circ(24)\\
& \mapsto\quad (12)\circ(14)\circ(34)\circ(13) 
&\mapsto\quad (12)\circ(23)\circ(12)\circ(24)\\
& \mapsto\quad (14)\circ(24)\circ(13)\circ(14) 
&\mapsto\quad (23)\circ(13)\circ(24)\circ(14)\\
& \mapsto\quad (14)\circ(13)\circ(24)\circ(14) 
&\mapsto\quad (23)\circ(24)\circ(13)\circ(14)\\
& \mapsto\quad (34)\circ(14)\circ(12)\circ(24) 
&\mapsto\quad (34)\circ(23)\circ(34)\circ(13)\\
& \mapsto\quad (34)\circ(12)\circ(24)\circ(24) 
&\mapsto\quad (34)\circ(34)\circ(24)\circ(13)
\end{align*}
The third and fourth factorizations belong to case (2) of the lemma.
These subfactorizations are changed in all four positions.
Every factor $(12)$ is replaced by $(34)$ and vice versa, and the same occurs
with $(13)$ and $(24)$.
But this amounts to the same result as  conjugation of each factor
by $\pi$.

The claim can be checked for the remaining 12 subfactorisations
by elementary computations of the same kind, or by
the observation that they correspond to the cases above under 
conjugation by $\pi$, by $(12)(34)$, or by $(13)(24)$.
\qed\\

To define a suitable invariant of a factorization in $\hat{\mathbb
O}$, we observe that
each factor of a  factorization $\tau$ in $\hat{\mathbb O}$ coincides
either with the
corresponding element of the  factorization $\tau^0$ or with its
conjugate by $\pi$.

The product of the  factorization is  left invariant by the braid group
action, whence  we claim that
$$
\#\{i \:|\: \tau_i\neq\tau_i^0\}/2 \quad \in \quad \N,
$$
and thus we get a well-defined parity invariant of the $\hat H$-action.

The claim is easily established in this way: changing one of the factors in the
right half amounts (in view of the commutativity of those factors) to
multiplication
by $(1,3)(2,4)$ on the right, while changing one of the factors in the
left half amounts  to multiplication
by $(1,2)(3,4)$ on the left. Since $(1,2)(3,4) z (1,3)(2,4) \neq z $
$\forall z \in \gS_4$ we conclude that the product is left
unchanged only if an even number of changes are made.

\begin{lemma}
The action of $\hat H$ on elements in $\hat{\mathbb O}$ preserves
the parity of the integer
$$
M \quad  : = \quad
\#\{i \:|\: \tau_i\neq\tau_i^0\}/2
\:\:+\:\:
\#\{i>4d,\:i\equiv 0 \ \rm{ ( mod \ 2)} \:|\: \tau_i\neq\tau_i^0\}
$$
\end{lemma}

\proof
All generators of $\hat H$ except $\sigma_s$ act either  trivially or
by transpositions
of factors on the same side (left or right) and at even distance,
    hence they preserve the first summand and likewise the second summand.

The action of $\sigma_s$ is different. We have the two possibilities of
lemma \ref{s-act}.
In the first case $\sigma_s$ acts trivially on $\tau$ and thus on $M$.

In the second case the product
$\tau_{4d-1}\tau_{4d}\tau_{4d+1}\tau_{4d+2}$ is equal to $(12)(34)$
or $(13)(24)$ which implies, since in the Klein group ($\cong (\Z/2)^2$) a product of two
non trivial elements is non trivial,  that this part of the factorization is one
of the following factorizations:
\begin{eqnarray*}
(12)\circ(34)\circ(13)\circ(13), &
(12)\circ(34)\circ(24)\circ(24),\\
(34)\circ(12)\circ(13)\circ(13), &
(34)\circ(12)\circ(24)\circ(24),
\end{eqnarray*}
resp.
\begin{eqnarray*}
(12)\circ(12)\circ(13)\circ(24),&
(12)\circ(12)\circ(24)\circ(13),\\
(34)\circ(34)\circ(13)\circ(24),&
(34)\circ(34)\circ(24)\circ(13).
\end{eqnarray*}
In all cases it can be checked that
$\#\{i\in\{4d-1,4d,4d+1,4d+2\} \:|\: \tau_i\neq\tau_i^0\}$
is either $1$ or $3$
and $\#\{i=4d+2 \:|\: \tau_i\neq\tau_i^0\}$ is either $0$ or $1$.

Both numbers change when the factorization is conjugated by $\pi$,
hence the parity of $M$ is preserved, as claimed.

\qed

Finally we can prove the essential step:

\begin{prop}\label{nonconj}
The elements $\sigma_p^2$ and $\sigma_q^2$ are not conjugate in
$\hat H$. Hence also $\sigma_p$ and $\sigma_q$ are not conjugate in
$\hat H$.
\end{prop}

\proof
Suppose on the contrary that there exists an $h\in \hat H$ such that
$\sigma_p h = h \sigma_q$.
Then the action of both sides on the  factorization $\tau^0$ must give
the same  factorization.
But the invariant $M$ is $2$ for $\tau^0h\sigma_q=\tau^0\sigma_q$
and odd for $\tau^0\sigma_p h$, since it is $1$ for $\tau_0\sigma_p$.
So the  factorizations are different and we get
a contradiction.

The second assertion was already shown at the beginning of the present section.

\qed

\section{Non equivalence of the horizontal Lefschetz fibrations in 
some cases. }

Consider an $abc$-surface $S$ and its vertical
  Lefschetz fibration $F:  S \ra \PP^1$ given by  the
composition of the perturbed covering $\pi : S \ra Q$
with the  first projection onto $ \PP^1$; the differentiable class of $F$ is encoded in the
M-equivalence class of the factorization of the identity in
the Mapping class group  $\MMM ap_g$   corresponding to the monodromy
factorization of $F$ ($g$ is  the genus of 
the base fibre $C_0$ of $F$) .

As already mentioned, we  proved in \cite{c-w} that if $S'$ is an 
$(a',b',c')$ surface, where $b=b'$ and $a+c=a'+c'$ then $S$ and $S'$ 
are
not deformation equivalent but the above defined Lefschetz fibrations 
are M-equivalent,
in particular proving that the surfaces $S$ and $S'$ are diffeomorphic via
an orientation preserving diffeomorphism sending the canonical class
to the canonical class.

We want to give here another small contribution in the direction
of the  interesting and difficult question whether the surfaces are 
symplectomorphic with respect to the canonical symplectic
structure induced by the complex structure.

The projections of $S$ and $S'$ to $Q = \PP^1\times\PP^1$ have 
discriminant curves $\Delta$ and $\Delta '$.
Our  computations (\cite{clw}) of the Braid monodromy factorizations 
of these curves with respect to the vertical
projection (onto the $x$-axis),
that we illustrated previously,
  show that the braid monodromy factorizations are not stably equivalent.
  This rises the suspicion that the surfaces may be not
canonically symplectomorphic.\par\medskip

We shall in fact point out another difference between the surfaces 
$S$ and $S'$.
We consider the Lefschetz fibrations corresponding to the horizontal 
projection (to the $y-axis$).
  This time a regular fibre $C_0'$ is a  4-tuple covering of  a
  horizontal line in $Q$, with  (4a+4c) branch points , hence it has genus
2a+2c-3.  It is identified in a natural way with the curve on Figure 
\ref{fibre}.

\unitlength=1pt
\begin{figure}
\centering
\begin{picture}(400,150)(-20,-40)

  \put(240,0){$\dots$}

   \put(110,0){$\dots$}
  \put(301,2){$.$}  \put(305,2){$.$} \put(309,2){$.$}
\put(311,2){$.$}  \put(315,2){$.$}
\put(61,2){$.$}  \put(65,2){$.$} \put(69,2){$.$}
\put(71,2){$.$}  \put(75,2){$.$}
\put(21,82){$.$}  \put(25,82){$.$} \put(29,82){$.$}
\put(31,82){$.$}  \put(35,82){$.$}
\put(261,82){$.$}  \put(265,82){$.$} \put(269,82){$.$}
\put(271,82){$.$}  \put(275,82){$.$}

        \put(220,80){$\dots$}  \put(90,80){$\dots$}
\put(31,80){\oval(22,10)[b]} \put(269,80){\oval(22,10)[b]}
\put(71,0){\oval(22,10)[b]}  \put(309,0){\oval(22,10)[b]}
  \put(150,80){\oval(260,40)[t]} \put(80,80){\oval(120,40)[bl]}
\put(80,40){\oval(120,40)[tr]} \put(110,40){\oval(60,40)[br]}
\put(110,0){\oval(100,40)[tl]} \put(190,0){\oval(260,40)[b]}
\put(260,0){\oval(120,40)[tr]} \put(260,40){\oval(120,40)[bl]}
\put(230,40){\oval(60,40)[tl]} \put(230,80){\oval(100,40)[br]}
\put(170,40){\oval(20,20)} 
\put(50,80){\oval(16,16)}
\put(200,80){\oval(16,16)} \put(140,80){\oval(16,16)}
  \put(90,0){\oval(16,16)} \put(140,0){\oval(16,16)}
\put(200,0){\oval(16,16)} \put(290,0){\oval(16,16)}
\put(250,80){\oval(16,16)}


\qbezier[5](143,72)(142,68)(150,58)
\qbezier[5](150,58)(158,48)(162,47)
\qbezier(145,73)(157,66)(164,49)

\qbezier[5](197,72)(198,68)(190,58)
\qbezier[5](190,58)(182,48)(178,47)
\qbezier(195,73)(183,66)(176,49)

\qbezier[5](197,8)(198,12)(190,22)
\qbezier[5](190,22)(182,32)(178,33)
\qbezier(195,7)(183,14)(176,31)

\qbezier[5](143,8)(142,12)(150,22)
\qbezier[5](150,22)(158,32)(162,33)
\qbezier(145,7)(157,14)(164,31)

\put(137,55){$\alpha_1$}
\put(193,55){$\beta_1$}
\put(190,25){$\gamma_1$}
\put(140,25){$\delta_1$}
\put(45.8,92){$\alpha_{2a-2}$} \put(135,92){$\alpha_2$} 
\put(242,92){$\beta_{2c-2}$}
\put(195,92){$\beta_2$} \put(168,36){$\sigma$} \put(195,-14){$\gamma_2$}
\put(282,12){$\gamma_{2a-2}$} \put(136,-16){$\delta_2$} 
\put(85,12){$\delta_{2c-2}$}
  \put(24,68){$\alpha_{2a-1}$} \put(254.5,67.8){$\beta_{2c\text{-}\!\;\!1}$}
\put(64.8,-12.8){$\delta_{2c-1}$} \put(294,-12){$\gamma_{2a\text{-}1}$}

\end{picture}
\caption{The fibre $C_0'$ of the horizontal Lefschetz fibration}\label{fibre}
\end{figure}
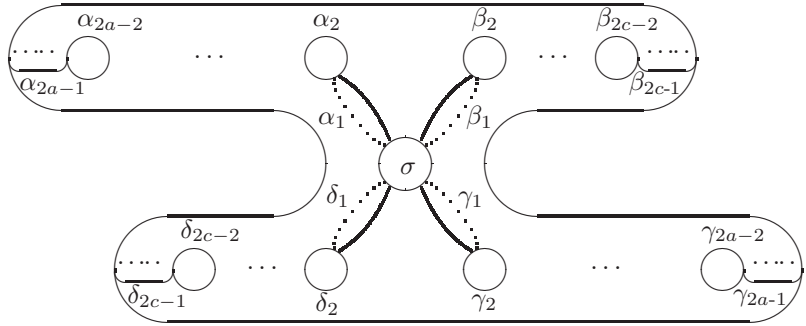

\unitlength=.53mm
\linethickness{.8pt}

\begin{figure}
\centering
\begin{picture}(240,100)(-140,-20)

\qbezier(-140,0)(-140,20)(-20,20)
\qbezier(100,0)(100,20)(20,20)
\qbezier(-140,0)(-140,-20)(-35,-20)
\qbezier(100,0)(100,-20)(35,-20)

\qbezier(-110,-5)(-105,-5)(-100,2)
\qbezier(-110,-5)(-115,-5)(-120,2)
\qbezier(-117,0)(-110,5)(-103,0)

\qbezier(-70,-5)(-65,-5)(-60,2)
\qbezier(-70,-5)(-75,-5)(-80,2)
\qbezier(-77,0)(-70,5)(-63,0)

\qbezier(70,-5)(65,-5)(60,2)
\qbezier(70,-5)(75,-5)(80,2)
\qbezier(77,0)(70,5)(63,0)

\qbezier(-3,3.5)(0,4)(3,3.5)
\qbezier(-38,0)(-30,3)(-18,3.5)
\qbezier(18,3.5)(30,3)(38,0)

\qbezier(-140,60)(-140,80)(0,80)
\qbezier(100,60)(100,80)(0,80)
\qbezier(-140,60)(-140,40)(-40,40)
\qbezier(100,60)(100,40)(40,40)

\qbezier(-110,55)(-105,55)(-100,62)
\qbezier(-110,55)(-115,55)(-120,62)
\qbezier(-117,60)(-110,65)(-103,60)

\qbezier(-70,55)(-65,55)(-60,62)
\qbezier(-70,55)(-75,55)(-80,62)
\qbezier(-77,60)(-70,65)(-63,60)

\qbezier(70,55)(65,55)(60,62)
\qbezier(70,55)(75,55)(80,62)
\qbezier(77,60)(70,65)(63,60)

\qbezier(-38,60)(0,65)(38,60)

\qbezier(40,40)(20,40)(0,8)
\qbezier(-35,-20)(-17.5,-20)(0,8)

\qbezier(-45,2)(-25,-8)(-20,0)
\qbezier(-20,0)(10,50)(45,62)

\put(8,19){\color{white}\circle*{5}}
\put(-8,19){\color{white}\circle*{5}}
\put(0,8){\color{white}\circle*{5}}
\put(0,29){\color{white}\circle*{5}}
\put(4,13.5){\color{white}\circle*{3}}
\put(-4,24){\color{white}\circle*{3}}

\qbezier(-40,40)(-20,40)(0,8)
\qbezier(35,-20)(17.5,-20)(0,8)

\qbezier(45,2)(25,-8)(20,0)
\qbezier(20,0)(-10,50)(-45,62)


\linethickness{.2pt}
\qbezier(-70,11)(-55,11)(-55,0)
\qbezier(-85,0)(-85,11)(-70,11)
\qbezier(-70,-12)(-55,-12)(-55,0)
\qbezier(-85,0)(-85,-12)(-70,-12)

\linethickness{.2pt}
\qbezier(-39,0)(-40,10)(-51,10)
\qbezier(-62,0)(-62,10)(-51,10)

\linethickness{1pt}
\qbezier[8](-39,0)(-40,-10)(-51,-10)
\qbezier[8](-62,0)(-62,-10)(-51,-10)

\linethickness{.2pt}
\qbezier(-118,0)(-118,8)(-126,7)
\qbezier(-140,0)(-140,5)(-126,7)

\linethickness{1pt}
\qbezier[5](-118,0)(-118,-8)(-126,-7)
\qbezier[8](-140,0)(-140,-5)(-126,-7)

\linethickness{2pt}
\qbezier[2](-89,0)(-93,0)(-97,0)

\put(-50,14){\makebox(0,0){$\gamma_1$}}
\put(-70,14.5){\makebox(0,0){$\gamma_2$}}
\put(-112,8){\makebox(0,0){$\gamma_{2a-1}$}}

\linethickness{.2pt}
\qbezier(-70,71)(-55,71)(-55,60)
\qbezier(-85,60)(-85,71)(-70,71)
\qbezier(-70,48)(-55,48)(-55,60)
\qbezier(-85,60)(-85,48)(-70,48)

\linethickness{.2pt}
\qbezier(-39,60)(-40,70)(-51,70)
\qbezier(-62,60)(-62,70)(-51,70)

\linethickness{1pt}
\qbezier[8](-39,60)(-40,50)(-51,50)
\qbezier[8](-62,60)(-62,50)(-51,50)

\linethickness{.2pt}
\qbezier(-118,60)(-118,68)(-126,67)
\qbezier(-140,60)(-140,65)(-126,67)

\linethickness{1pt}
\qbezier[5](-118,60)(-118,52)(-126,53)
\qbezier[8](-140,60)(-140,55)(-126,53)

\linethickness{2pt}
\qbezier[2](-89,60)(-93,60)(-97,60)

\put(-50,74){\makebox(0,0){$\alpha_1$}}
\put(-70,74.5){\makebox(0,0){$\alpha_2$}}
\put(-112,68){\makebox(0,0){$\alpha_{2a-1}$}}

\linethickness{.2pt}
\qbezier(39,60)(40,70)(51,70)
\qbezier(62,60)(62,70)(51,70)

\linethickness{1pt}
\qbezier[8](39,60)(40,50)(51,50)
\qbezier[8](62,60)(62,50)(51,50)

\linethickness{2pt}
\qbezier[2](84,0)(88,0)(92,0)

\put(50,74){\makebox(0,0){$\beta_1$}}

\linethickness{.2pt}
\qbezier(39,0)(40,10)(51,10)
\qbezier(62,0)(62,10)(51,10)

\linethickness{1pt}
\qbezier[8](39,0)(40,-10)(51,-10)
\qbezier[8](62,0)(62,-10)(51,-10)

\linethickness{2pt}
\qbezier[2](84,60)(88,60)(92,60)

\put(50,14){\makebox(0,0){$\delta_1$}}

\linethickness{.2pt}

\qbezier(-45,65)(0,75)(45,65)
\qbezier(-45,65)(-58.5,62)(-45,56)
\qbezier(45,65)(58.5,62)(45,56)
\qbezier(-15,-3)(-21,-13)(-45,-2)
\qbezier(-45,-2)(-49,-0.3)(-48,4)
\qbezier(-48,4)(-47,9)(-15,10)

\qbezier(15,-3)(-10,40)(-45,56)
\qbezier(15,-3)(21,-13)(45,-2)
\qbezier(45,-2)(49,-0.3)(48,4)
\qbezier(48,4)(47,9)(15,10)

\qbezier(-15,-3)(-10,6)(-4,14)
\qbezier(4,24)(28,49)(45,56)

\put(0,74){\makebox(0,0){$\sigma$}}

\end{picture}
\caption{}

\end{figure}

Our computation of the Braid monodromy of the  discriminant curve, 
where the roles of the
variables
$x$ are easily exchanged,
  yields easily the monodromy factorization
corresponding to the horizontal  Lefschetz fibration. In particular 
it shows that the monodromy group of the Lefschetz
fibration is generated by Dehn twists with respect to the cycles 
$\alpha_i,\beta_i,\gamma_i,\delta_i,\sigma$ on Figure 1.
  We shall prove that if $a+c=a'+c'$ is even and $a-a'$ is odd and 
$b=b'$ then the monodromy factorizations
   for the horizontal Lefschetz fibrations are not M-equivalent 
although the fibrations have the same genus and the surfaces $S$ and 
$S'$ are diffeomorphic.

We map the cycles in $C_0'$ into the homology group $H_1(C_0',\Z)$ and 
further into  the vector space $V=H_1(C_0',\Z / 2\Z)$ of vectors modulo 
$2$.
  The intersection form of cycles on $C_0'$ induces an alternating 
symmetric form $(u,v)$ on $V$. The mapping class group
   of $C_0'$ projects onto the symplectic group of $V$ with respect to 
the form and Dehn twists project to transvections $T$ with respect to 
cycles:
    $T_u(v)=v+(u,v)u$.

We shall denote the images of 
$\alpha_i,\beta_i,\gamma_i,\delta_i,\sigma$  in $V$ by 
$a_i,b_i,c_i,d_i,s$ respectively.

Considering the intersection matrix of these vectors we see that if 
we omit the vectors $b_{2c-1},c_{2a-1},d_{2c-1}$
  then the remaining vectors have a non-singular intersection matrix 
and thus they form a basis of $V$.

  We define a quadratic function $q:V\to \Z / 2\Z$ using this basis :

   $q$ is equal 1 on each basis vector and $q(u+v)=q(u)+q(v)+(u,v)$ 
for any pair of vectors $u,v$.

  Quadratic functions with an associated bilinear form as above fall into 
two equivalence classes according to their Arf invariant (modulo 2).
In order to compute the Arf invariant we choose a symplectic basis 
$e_i,f_i$ of $V$, such that $(e_i,f_i)=1$ for all $i$,
while all other intersections products of basis vectors are equal 0. 
Then the Arf invariant of $q$ is equal to $\Sigma_i  q(e_i)q(f_i)$.

  We now recall Theorem 4 from \cite{waj1}:

\begin{teo} Let $G$ be a subgroup of $Sp(V,\Z / 2\Z)$ generated by 
transvections, containing the transvections with respect to the 
vectors  $u_i$
of a fixed basis.
  Let $q$ be the quadratic function defined by $q(u_i)=1$.

If $G$ contains a transvection with respect to a vector $v$ 
satisfying $q(v)=0$ then $G=Sp(V,\Z / 2\Z)$.

If this is not the case and the basis is not special, e.g. the 
intersection diagram of the basis is a tree
but it is different from the Dynkin diagrams $A_n$ and $D_n$, then 
$G$ is equal to the orthogonal group of the function $q$,
the group of all elements of $Sp(V,\Z / 2\Z)$ preserving $q$.
\end{teo}
\medskip
Let us go back to the Lefschetz fibration of the surface $S$. Let $G$ 
be the image of the monodromy group, the group
generated by the transvections with respect to  the vectors 
$a_i,b_i,c_i,d_i,s$. Consider the basis described before.
The intersection diagram of this basis has the form of a cross, thus the basis is not special. 

Consider the corresponding quadratic function $q$.
The  vector
$b_{2c-1}$, which is not in the basis, is equal to the sum 
$v=a_1+a_3+ a_5 + \dots + a_{2a-1}+b_1+b_3+b_5 + \dots + b_{2c-3}$
(both vectors have the same intersection with each basis vector so 
they are equal).

 If $a+c$ is odd then $q(v)=0$
and the group $G$ generated by the transvections coincides with 
$Sp(V,\Z / 2\Z)$.

If $a+c$ is even then $q(b_{2c-1})=q(c_{2a-1})=q(d_{2c-1})=1$ and the 
group $G$ is equal to the orthogonal
group of $q$. In this case if we start with any basis $v_i$ 
satisfying $q(v_i)=1$
  and we define a quadratic function using the basis $v_i$ we get the 
same function $q$.

It is easy to compute the Arf invariant of $q$.

  To simplify notation we work out  an example with $a=2,c=2$.

We can choose the following symplectic basis : 
$(a_3,a_2),(a_3+a_1,s),(a_3+a_1+b_1,b_2),(a_3+a_1+c_1,c_2),(a_3+a_1+d_1,d_2)$. 

The sum of $k$ pairwise non-intersecting basis vectors has the 
$q$-value 1 if $k$ is odd and has the $q$-value 0 if $k$ is even.

Hence  the Arf invariant of $q$ is 0 if $a$ is even and is 1 if $a$ is odd.

We now prove that if $a+c=a'+c'$ is even and $a-a'$ is odd then the 
monodromy factorizations
  of the horizontal Lefschetz fibrations of an abc-surface $S$ and of an a'b'c'-surface $S'$ are not
  M-equivalent. Without loss of generality assume that a is even, and that a' is odd: so that the Arf invariant of $q$
  is $0$, while the Arf invariant of $q'$ is $1$.
  
  Assume the contrary: then, in fact,  making an  arbitrary identification of the 
fibers of the fibrations, we may  assume
   that we have two factorizations by Dehn twists which are Hurwitz equivalent.

    We map the corresponding cycles to vectors in $V$.
    The monodromy factorization maps to a factorization in $Sp(V,\Z / 2\Z)$ by 
transvections with respect to some vectors.

    This set of vectors contains a basis with the intersection pattern 
which defines a quadratic
    function $q$ with the Arf invariant  0. All vectors of the initial 
set have the $q$-value 1.
The action induced by a Hurwitz move in $Sp(V,\Z / 2\Z)$ replaces a 
pair of transvections
$T_u,T_v$ by the pair $T_v,T_v^{-1}T_uT_v$. One checks that 
$T_v^{-1}T_uT_v=T_w$
where $w=u+(u,v)v$. Assume $q(u)=q(v)=1$: then $q(w)=1$. So, after 
the Hurwitz
move, we have again a product of transvections with respect to 
vectors of $q$-value 1.
  If we assume that the factorizations are M-equivalent, then we get  a factorization 
by  transvections with
   respect to a set of vectors having $q$-value 1 and which contains a 
basis with an intersection pattern such that it
     defines a quadratic function $q'$ with Arf invariant 1. 
But, by the previous remark , $q'=q$. This is a contradiction.

One can easily show that the orthogonal groups of $q$ and $q'$ are 
not conjugate in $Sp(V,\Z / 2\Z)$
and therefore the monodromy groups of the two fibrations are not 
conjugate in the mapping class group of a fibre.
Probably these monodromy groups  are not isomorphic at all, but we do 
not know at the moment how to prove it.

\section{ Epilogue?}

The question of determining whether
the Braid Monodromy factorizations of two diffeomorphic surfaces $S$, 
$S'$ (an abc-surface, resp. an a'b'c'-surface)
are M-equivalent  was motivated by the related  question of
  deciding about the existence of
a diffeomorphism $\varphi$ between $S, S'$  commuting with the 
respective vertical Lefschetz fibrations,
and yielding a canonical
symplectomorphism.

Our main result that these factorizations are not stably equivalent 
suggests that such a diffeomorphism cannot exist.

  The relation of stable equivalence (the authors call it  indeed m-equivalence)
was introduced by Auroux
and Katzarkov in \cite{a-k1} in order to obtain invariants of
symplectic 4-manifolds.

The reason to introduce it,
allowing not only Hurwitz equivalence and simultaneous conjugation, but also
creation/cancellation of admissible
pairs of a positive and of a negative node (here a node, and
its   corresponding full twist $\beta$, is said
to be admissible if the inverse image of the node inside the
ramification divisor consists of two disjoint smooth branches)
has to do with the fact that  this phenomenon occurs
when considering isotopies of approximately holomorphic maps.

Our  result (\ref{bmf}) shows  that the braid monodromy group $H$ depends
only upon the numbers $ b $ and $d$, provided for instance that we have
nonvanishing  of the  respective numbers
$8ab-2(ad+bc)$,
$8cd-2(ad+bc)$,
or provided that we are in the case $b=d$.

A fortiori if the groups $H$ are the same for different choices of $(a,b,c,d)$
then the  fundamental groups $\pi_1 ( Q \setminus \De) $ are isomorphic.

Our new method  of distinguishing factorizations, and not only
associated fundamental groups of complements of branch curves (by the 
Zariski-van Kampen method),
leads to  the above results representing  the first
positive step towards the realization
of a more general program set up by Moishezon (\cite{mois1}, \cite{mois2})
in order to produce braid monodromy invariants which should distinguish
the connected components of a moduli space $\MM_{\chi,K^2}$.

Moishezon's program is based on the consideration (assume here for simplicity
that $K_S$ is ample) of
a general projection $\Psi_m : S \ra \PP^2$ of a pluricanonical
embedding $\Phi_m : S \ra \PP^{P_m -1}$, and of the braid monodromy
factorization
corresponding to the (cuspidal) branch curve $B_m$ of $\Psi_m$.

An invariant of the connected component is here given by the M-equivalence class
(i.e., for Hurwitz equivalence plus simultaneous conjugation) of this
braid monodromy factorization. Moishezon, and later Moishezon-Teicher,
calculated a coarser invariant, namely the fundamental group
$\pi_1(\PP^2- B_m)$.

This group turned out to be not overly complicated,
and in fact, as shown in many cases in \cite{a-d-k-y}, it tends  to give
no extra information beyond the one given by the topological
           invariants of $S$ (such as $\chi, K^2$).

Auroux and Katzarkov showed that, for $ m >>0$, the stable-equivalence class
of the above braid monodromy factorization determines the canonical
symplectomorphism class of $S$, and conversely.

The work by Auroux, Katzarkov adapted Donaldson's techniques for proving
the existence of symplectic Lefschetz fibrations (\cite{don6},\cite{donsympl})
in order to show  that each integral symplectic $4$-manifold is in a
natural way 'asymptotically' realized by  a generic symplectic covering of
$\PP^2_{\CC}$, given by approximately holomorphic sections of a high multiple
$L^{\otimes m}$ of a complex line bundle $L$
whose class is the one of the given
integral symplectic form.

The methods of  Donaldson on one side, Auroux and Katzarkov on the other,
     use algebro geometric methods in order to produce invariants of
symplectic 4-manifolds.

For instance, in the case  of a generic symplectic covering of the
plane, we get a corresponding branch curve $\De_m$  which is a symplectic
submanifold with singularities only nodes and cusps.

To  $\De_m$ corresponds then a factorization in the braid group,
called m-th braid monodromy factorization: it  contains, as in the 
present paper,
only factors which are conjugates of $ \sigma_1^j$, not only with $ j
=  1, 2,3$
as in the complex algebraic case,
but also with $ j = -2$.

The factorization is not unique because it may happen that a pair of
two consecutive nodes, one positive and one negative, may  be created,
or disappear, and this is the reason to consider its  stable 
equivalence class which
      the authors show , for $ m >> 1$, to be  an invariant
of the integral symplectic manifold.

As we have already remarked, in the case of abc-surfaces
the vertical Braid monodromy groups $H$ are determined
by $b$ (up to conjugation), hence the fundamental groups  $\pi_1 ( Q
\setminus \De) $ are isomorphic for a fixed value of $b$.

So, Moishezon's
technique produces no invariants.

Our result indicates that one should try to go all the way
in the direction of  understanding stable-equivalence classes of pluricanonical
braid monodromy factorizations.  Let us describe here how
this could be done.

Define  $ p : S \ra \PP^2$ as the morphism  $p$
obtained as the composition of
$\pi : S \ra Q$ with the embedding $ Q \hookrightarrow \PP^3$
followed by a general projection  $\PP^3 \dashrightarrow \PP^2$.

In the even more special case of abc-surfaces such that $ a+c = 2b$,
           the m-th pluricanonical mapping $\Phi_m : S \ra \PP^{P_m -1}$
has a (non generic) projection given by the composition of $p$ with
a Fermat type map  $\nu_r : \PP^2  \ra \PP^2$
(given by $ \nu_r (x_0, x_1, x_2) = (x_0^r, x_1^r, x_2^r) $
in a suitable linear coordinate system),
where $ r : = m ( 2b -2) $.

Let $B$ be the branch curve of a generic perturbation of $ p$:
then the braid monodromy factorization
corresponding to $B$ can be calculated from the braid monodromy factorization
corresponding to $\De$.

We hope, in a sequel to this paper,
to be able to determine whether these braid monodromy factorizations are
equivalent, respectively stably-equivalent,  for abc-surfaces such that
$ a+c = 2b$.

The problem of calculating the braid monodromy factorization corresponding
to the (cuspidal) branch curve $B_m$ starting from the
braid monodromy factorization of $B$ has been addressed, in the
special case $ m=2$,
by Auroux and Katzarkov (\cite{a-k2}). Iteration of
their formulae should lead to the calculation of
           the braid monodromy factorization corresponding
to the (cuspidal) branch curve $B_m$ in the case, sufficient for applications,
           where $m$ is a sufficiently large power of $2$.

 Whether it will be possible to decide about the 
stable equivalence  of the factorizations  which  these formidable calculations will yield,
it is a quite
open question.  

In both directions the result would be extremely
interesting, leading either to

i) a counterexample to the speculation DEF = CAN. SYMPL also in the simply
connected case, or to

ii) examples of diffeomorphic but not canonically symplectomorphic
simply connected algebraic surfaces.

{\bf Acknowledgements:} 
The present research was performed in the realm of the Forschergruppe 790
`Classification of algebraic surfaces and compact complex manifolds'
of the DFG (DeutscheForschungsGemeinschaft).

We are grateful to a referee for a very careful reading of the manuscript which saved us
from  many  misprints and inaccuracies.

\smallskip


\end{document}